\documentclass[twoside]{article}
\usepackage{imsart}
\usepackage{amsfonts,amsmath,latexsym,amsthm,amssymb}

\newtheorem{prop}{Proposition}[section]
\newtheorem{rema}[prop]{Remark}

\newtheorem{defi}[prop]{Definition}
\newtheorem{assu}[prop]{Assumption}
\newtheorem{lemm}[prop]{Lemma}
\newtheorem{coro}[prop]{Corollary}
\newtheorem{exam}[prop]{Example}
\newcommand{\filtrationf}{$\mathbb{F}$}
\newcommand{\filtrationg}{$\mathbb{G}$}
\newcommand{\beqn}{\begin{eqnarray}}
\newcommand{\eeqn}{\end{eqnarray}}
\newcommand{\nonu}{\nonumber}

\newcommand{\pqa}{\left[}
\newcommand{\pqc}{\right]}
\newcommand{\pga}{\left\{}
\newcommand{\pgc}{\right\}}
\newcommand{\pta}{\left(}
\newcommand{\ptc}{\right)}
\newcommand{\vaa}{\left|}
\newcommand{\vac}{\right|}
\newcommand{\shl}{\mathcal{L}}
\newcommand{\shg}{\mathcal{G}}
\newcommand{\shf}{\mathcal{F}}
\newcommand{\shd}{\mathcal{D}}
\newcommand{\sha}{\mathcal{A}}
\newcommand{\shm}{\mathcal{M}}
\newcommand{\ep}{\varepsilon}
\newcommand{\interval}{0\leq t \leq 1}
\newcommand{\intervala}{0\leq t< 1}
\begin{document}

\begin{frontmatter}
\title{Modeling financial assets without semimartingale}

\author{Rosanna COVIELLO,$^{a,b}$}
\author{Francesco RUSSO$^{a}$}
\address{a Universit\'{e} Paris 13, Institut Galil\'{e}e,
Math\'ematiques, 99, avenue J.B.~Cl\'{e}ment, F-93430
Villetaneuse, France}
\address{b Modal'X, Universit\'{e} Paris 10, Avenue de la R\'epublique 92001 Nanterre Cedex, France}
\runauthor{Rosanna Coviello and Francesco Russo}

\begin{center}
\address{June 26, 2006}
\end{center}

\hskip2cm

\begin{abstract}
This paper does not suppose a priori that the evolution of the price of a financial asset
is a semimartingale. Since possible strategies of investors are self-financing,
previous prices are forced to be finite quadratic variation processes.
The non-arbitrage property is not excluded if the class $\sha$ of admissible strategies
is restricted. The classical notion of martingale is replaced with the notion
of $\sha$-martingale. A calculus related to  $\sha$-martingales with some examples
is developed. Some applications to the maximization of the utility of an
insider are expanded.
\end{abstract}


\bigskip


\end{frontmatter}

\textwidth15cm
\oddsidemargin1cm
\parindent0mm
\parskip1ex
\thispagestyle{empty}
\markright{}

\oddsidemargin =10mm \evensidemargin =10mm \topmargin =5mm
\textwidth =150mm \textheight =200mm

\section{Introduction}
{According to 
the fundamental theorem of asset pricing of Delbaen and Schachermayer in \cite{DS},
 in absence of \textit{free lunches with vanishing risk} (NFLVR), when investing possibilities run only trough simple predictable strategies  with respect to some filtration $\mathbb{G}$,
the price process of the risky asset $S$  is forced to be a semimartingale.
However (NFLVR) condition could not be reasonable in several situations. In that case $S$ may not be a semimartingale.
We illustrate here some of those circumstances.

Generally, admissible strategies are let vary in a quite large class of
predictable processes with respect to some filtration $\mathbb{G}$, representing the information flow available to the investor. 
As a matter of fact, the class of admissible strategies could be reduced because of  different market regulations 
or for practical reasons. For instance, the investor  could not be allowed to hold more than a certain number of stock shares.
On the other hand it could be realistic to impose a minimal delay between
two possible transactions as suggested by Cheridito (\cite{che}): when the logarithmic price $\log(S)$ is a geometric fractional Brownian motion (fbm), it is impossible to realize arbitrage possibilities satisfying that minimal requirement. We remind that without that restriction,
the market admits arbitrages, see for instance \cite{rogers}.
When the logarithmic price of $S$ is a geometric fbm or some particular strong Markov process, arbitrages can be excluded taking into account proportional transactions costs: Guasoni (\cite{gua}) has shown that, in that case,  the class of admissible strategies has to be restricted to bounded variation processes and this rules out arbitrages.

Besides the restriction of the class of admissible strategies, the adoption of non-semimartingale models finds its justification when the no-arbitrage condition itself is not likely. 

Empirical observations  reveal, indeed, that $S$ could fail to be a semimartingale because of market imperfections due to micro-structure noise,
as intra-day effects.
A model which considers those imperfections 
 would add to  $W$, the Brownian motion describing  $\log$-prices, a zero quadratic variation process, as a fractional Brownian motion of Hurst index greater than $ \frac{1}{2}$,
see for instance \cite{wor}.
Theoretically arbitrages in very small time interval could be possible, which would be compatible
with the lack of semimartingale property.

At the same way if (FLVR) are not possible for an \textit{honest investor}, an {\it inside trader}
could realize a free lunch with respect to the  enlarged filtration $\mathbb{G}$ 
including the one generated by prices and the extra-information.
Again in that case  $S$ may not be a semimartingale.
The literature concerning inside trading and asymmetry of information has been extensively enriched by several papers in the last ten years;
among them we quote Pikowski and Karatzas (\cite{pikar}),  Grorud and Pontier (\cite{gropon}),
 Amendinger, Imkeller and Schweizer (\cite{AIS}). They adopt enlargement of filtration techniques to describe the  evolution of stock prices in the insider filtration. 

Recently, 
some authors approached the problem in a new way using in particular forward integrals, in the framework of 
stochastic calculus via regularizations. For a comprehensive survey of that calculus  see \cite{RV05}.
Indeed, forward integrals could  exist also for
non-semimartingale integrators.
Leon, Navarro and Nualart in \cite{LNN}, for instance, solve the problem of maximization of expected logarithmic utility of an agent who holds an initial information depending on the future of prices.
They operate under  technical conditions which, a priori, do not imply the classical  assumption (H') for enlargement considered 
in \cite{jacod}. 
Using  forward integrals, they determine the utility maximum.
However, a posteriori, they found out that  their conditions let $\log(S)$ be a semimartingale.

Biagini and {\O}ksendal (\cite{BO}) considered somehow the converse implication. Supposing that the maximum utility is attained, they proved that
$S$ is a semimartingale.
Ankkirchner and Imkeller (\cite{AI})  continue to develop the enlargement of filtrations techniques and show, among the others, a similar result as \cite{BO}
using the fundamental theorem of asset pricing  of Delbaen-Schachermayer. In particular they establish a link 
between that fundamental theorem 
and finite utility. 

In our paper we treat a market where there are one risky asset, whose price is a strictly positive process $S$,
and a \textit{less risky} asset with price $S^0$, possibly  riskless but a priori only with bounded variation. A class $\sha$ of admissible trading strategies is 
specified. If $\sha$ is not large enough to generate all predictable simple strategies, then $S$ has no need to be a semimartingale,
even requiring the absence of free lunches among those strategies.
We try to build the basis of a corresponding financial theory which allows to deal with several problems as hedging  and non-arbitrage pricing, viability and completeness as well as with utility maximization. 

For the sake of simplicity in this introduction we suppose that the less risky  asset $S^0$ is constant and equal to 1.

As anticipated, a natural tool to describe the self-financing condition
is the forward integral of an integrand process $Y$ with respect to an integrator $X$,  denoted by $\int_0^t Y d^-X$; see section \ref{s2} for definitions.
Let $\mathbb{G} = (\shg_t)_{0 \le t \le 1}$ be a filtration on an underlying probability space $(\Omega, \shf, P)$, with $\shf = \shg_1$; $\mathbb{G}$
represents the flow of information available to the investor.
A {\bf self-financing portfolio} is a pair $(X_0,h)$ where $X_0$ is the initial value of the portfolio and $h$ is a  $\mathbb{G}$-adapted   and
 $S$-forward integrable process  specifying 
the number of  shares of  $S$ held in the portfolio. The market value process $X$ of such a portfolio, is given by $X_0+\int_0^\cdot h_sd^-S_s$, while 
$h^0_t = X_t - S_t h_t$ constitutes the number of shares of the less risky asset held.

 This  formulation of  self-financing condition is coherent with the discrete-time case. Indeed, let we consider a \textit{buy-and-hold strategy},
 i.e. a pair $(X_0,h)$ with $h= \eta I_{(t_0,t_1]}, 0 \le t_0 \le t_1 \le 1, $  and $\eta$ being a $\shg_{t_0}$-measurable random variable.
 Using the definition of forward integral it is not difficult to see that:
 $X_{t_0}=X_0,$  $X_{t_1}=X_0+\eta(S_{t_1}-S_{t_0})$. This implies   $h^0_{t_0^+}=X_0-\eta S_{t_0},$ $h^0_{t_1^+}=X_0+\eta(S_{t_1}-S_{t_0})$ and   
\beqn \label{1}
X_{t_0}=h_{{t_0}^{+}}
S_{t_0}+h^0_{t_0^+},\quad X_{t_1}=h_{t_1^+}S_{t_1}+h^0_{t_1^+}:
\eeqn
at the \textit{re-balancing} dates $t_0$ and $t_1,$ the value of the old portfolio must be reinvested to build the new portfolio without exogenous withdrawal of money.
    
In this paper $\sha$ will be a real linear subspace of all  self-financing portfolios  and it will constitute, by definition, the class of all {\bf admissible} portfolios.
 $\sha$  will depend on the kind of problems one has to face:
hedging, utility maximization,  modeling inside trading.
If we require that $S$ belongs to $\sha$, then the process $S$ is forced to be a finite quadratic variation process.
In fact, $\int_0^\cdot S d^-S$ exists if and  only if the quadratic variation $[S]$ exists, see \cite{RV05}; in particular
one would have 
$$  
\int_0^\cdot S_s d^- S_s = S^2 - S^2_0 - \frac{1}{2} [S]. 
$$

$\shl$ will be the sub-linear space of $L^0(\Omega)$ 
representing a  set of \textbf{contingent claims of interest} for one investor.
An $\sha$-{\bf attainable contingent claim} will be  a random variable $C$ for which  there is a self-financing portfolio  $(X_0, h)$ with $h \in \sha$ and 
$$ C = X_0 + \int_0^1 h_s d^- S_s. $$ $X_0$ will be called {\bf replication price} for $C$.
The market will be said $(\mathcal{A},\shl)$-{\bf complete} if every element of $\shl$ is $\sha$-attainable.

In these introductory lines we will focus only on
one particular elementary situation.

For simplicity we illustrate the case where $[\log (S)]_t=\sigma^2 t. $    
We choose as $\shl$ the set of all \textit{European} contingent claims $C = \psi(S_1) $ where $\psi$ is continuous with polynomial growth. 
We consider the case $\sha=\mathcal{A}_S,$ where 
\beqn \nonu
\mathcal{A}_S&=&\pga (u(t,S_t)), \intervala
\left |\right. u:[0,1]\times \mathbb{R} \rightarrow \mathbb{R}, \mbox{ Borel-measurable }\right. \\ \nonu &&\left.
\mbox{ with polynomial growth and lower bounded} \pgc.
\eeqn 
Such a market is $(\mathcal{A},\shl)$-complete: in fact, a random variable  $C = \psi(S_1)$ is an $\sha$-attainable contingent claim. To build a replicating strategy the investor has to 
choose $v$ as solution of the following problem 
$$ \left \{
\begin{array}{lll}
\partial_t v(t,x) + \frac{1}{2}\sigma^2 x^2  \partial^{(2)}_{xx} v(t,x) &= & 0 \\
v(1,x) &=& \psi(x) 
\end{array}
\right. $$
and $X_0 = v(0,S_0).$
This follows easily after application of It\^o formula contained in proposition \ref{p2.1}, see proposition \ref{p5.31}.

We highlight that this method can be adjusted to hedge also \textit{Asian} contingent claims.

A crucial concept is the one of
 $\sha$-martingale processes. Those processes naturally intervene in utility maximization, arbitrage and 
uniqueness of hedging prices.

A process $M$  is said to be an  $\sha$-{\bf martingale}
if for any process $Y \in \sha$, 
$$
E \pqa \int_0^\cdot Y d^-M \pqc= 0.
$$
If for some filtration   $\mathbb{F}$ with respect to which $M$ is adapted, 
   $\sha$ contains the class of all bounded $\mathbb{F}$-predictable processes, then $M$ is an $\mathbb{F} $-martingale.

An example of $\sha$-martingale is the so called {\bf weak Brownian motion of order} $k = 1$ and quadratic variation equal to $t$.
That notion was introduced in \cite{FWY}: a weak Brownian motion of order $1$ is a process $X$ such that the law of $X_t$ is $N(0,t)$ for any $t \ge 0$.

A portfolio $(X_0,h)$ is said to be an $\sha$-arbitrage if $h\in \sha$, $X_1 \ge X_0 $ almost surely and
$ P\{X_1 - X_0 > 0\} > 0$. 
We denote by $\shm$ the set of probability measures being equivalent to the initial probability  $P$  under which 
$S$ is an $\sha$-martingale.
If $\shm$ is non empty then the market is $\sha$-arbitrage free. 
In fact if  $Q \in \shm$, 
 given a pair $(X_0,h)$ which is an $\sha$-arbitrage, then $E^Q [X_1 - X_0] = E^Q [\int_0^1 h d^-S] = 0$.
In that case the replication price $X_0$ of an $\sha$-attainable contingent claim $C$ is unique,
 provided that the process $h\eta,$
for any bounded random variable $\eta$ in $\shg_0$ and $h$ in $\cal A,$ still  belongs to $\sha$.
Moreover  $X_0 = E^Q [C \vert \shg_0]$.
In reality, under the weaker assumption that  the market is $\sha$-arbitrage free,  the replication price is still unique, see proposition \ref{p5.33}.
Furthermore if $\mathcal{M}$ is non empty and $\sha = \sha_S,$ as assumed in this section, the law of $S_t$ has to be equivalent to Lebesgue measure for every $\interval,$ see 
proposition  \ref{p5.27}.

If the market is $(\mathcal{A},\shl)$-complete then all the probabilities measures in $\shm$ coincide on $\sigma(\shl)$, 
see proposition \ref{p5.30}.
If  $\sigma(\shl) = \shf$ then $\shm$ is a singleton: this result recovers the classical case.

Given an utility function satisfying usual assumptions, it is possible to show that the
maximum $\pi$ is attained on a class of portfolios  fulfilling
conditions related to assumption \ref{a5.40},
if and only if there exists a probability measure under which $\log(S)-\int_0^\cdot \pta \sigma^2 \pi_t-\frac{1}{2}\sigma^2\ptc dt$ is an $\cal{A}$-martingale, see proposition \ref{p5.45}. Therefore if $\sha$ is big enough to fulfill
conditions related to assumption $\shd$ in Definition \ref{d3.6}, 
then $S$ is a classical semimartingale.

Those considerations show that
most of the classical results of basic financial  theory admit a natural extension to
 non-semimartingale models.

The paper is organized as follows. After some preliminaries about stochastic calculus 
via regularizations for forward integrals,
we provide in section 3 examples of integrators and integrands for which forward integrals exist
and realize some important properties in view of financial applications:
those examples appear in three essential situations coming from Malliavin calculus, substitution formulae
and It\^o-fields. Regarding finance applications, the class of strategies defined using Malliavin calculus 
are useful when   $\log(S)$ is a geometric Brownian motion with respect to a filtration $\mathbb{F}$ contained
in $\mathbb{G}$; the use of substitution formulae naturally appear when trading 
with an initial extra information, already available at time $0$; 
 It\^o fields apply whenever $S$ is a generic finite quadratic variation process.

Section 4 is devoted to the study of $\sha$-martingales: after having defined and established basic properties,
we explore the relation between $\sha$-martingales and weak Brownian motion; later we discuss
the link between the existence of a maximum for a an optimization problem and the $\sha$-martingale property.

In Section 5 we finally deal with applications to mathematical finance. 
We define self-financing portfolio strategies and we provide examples.
Moreover we face technical problems related to the use of forward integral in order to 
describe the evolution of the wealth process. Those problems arise because of the lack of
chain rule properties.
Later, we discuss absence of $\sha$-arbitrages, $(\mathcal{A},\shl)$-completeness and hedging.
We conclude the section analyzing the problem of maximizing expected utility from terminal wealth. We obtain results about the existence of an \textit{optimal portfolio} generalizing those of \cite{LNN} and \cite{BO}.}

\section{Preliminaries}\label{s2}
For the convenience of the reader we give some basic concepts and
fundamental results about stochastic calculus with respect to finite
quadratic variation processes which will be extensively used later.
For more details we refer the reader to \cite{RV05}.

In the whole paper $\pta \Omega,
\mathcal{F},P\ptc$ will be a fixed probability space.
For a stochastic process $X=(X_t,\interval)$  defined on $\pta \Omega,
\mathcal{F},P\ptc$ we will adopt the convention $X_t=X_{(t\vee
  0)\wedge 1},$ for $t$ in $\mathbb{R}.$
Let $0\leq T\leq 1.$ 
We will say that a sequence of processes $\pta X^n_t, 0\leq t\leq T \ptc_{n\in_{\mathbb{N}}}$ \textbf{converges uniformly in probability (ucp) on $[0,T]$} toward a process $(X_t, 0\leq t\leq T),$ if $\sup_{t\in[0,T]}\vaa X^n_t-X_t\vac$ converges to zero in probability. 

\begin{defi}
\begin{enumerate}
\item
Let $X=(X_t,0\leq t\leq T)$ and $Y=(Y_t,0\leq t\le T)$ be processes with paths respectively in $C^0([0,T])$ and $L^1([0,T])$. Set, for every $0\leq t \leq T$,
$$
I(\ep,Y,X,t)=\frac{1}{\ep}\int_0^t Y_s \pta X_{s+\ep}-X_s \ptc ds, 
$$
and 
$$
C(\ep,X,Y,t)= \frac{1}{\ep} \int_0^t\pta Y_{s+\ep}-Y_s\ptc \pta X_{s+\ep}-X_s\ptc ds.
$$
If $I(\ep,Y,X,t)$ converges in probability for every $t$ in $[0,T],$ and the limiting process admits a continuous version $I(Y,X,t)$ on $[0,T],$
${Y}$ is said to be $\textbf{X}$\textbf{-forward integrable} on $[0,T]$. The process $\pta I(Y,X,t), 0\leq t\leq T\ptc$ is denoted by $\int_0^\cdot Yd^-X.$ 
If $I(\ep,Y,X,\cdot)$ converges $ucp$ on $[0,T]$ we will say that the forward integral $\int_0^\cdot Yd^-X$ is the \textbf{limit ucp of its regularizations}.
\item 
If $(C(\ep,X,Y,t), 0\leq t\leq T)$ converges ucp on $[0,T]$ when $\ep$ tends to zero, the limit will be called the \textbf{covariation process} between $X$ and $Y$   and it will be denoted by  $[X,Y].$
If $X=Y,$ $\pqa X,X\pqc$ is called the \textbf{finite quadratic variation} of $X$: it will also be denoted by $\pqa X\pqc,$ and $X$ will be said to be a \textbf{finite quadratic variation process} on $[0,T].$ 
\end{enumerate}  
\end{defi}

\begin{defi}\label{d3.2} We will say that a process  $X=(X_t,0\leq t\leq T),$ is
  \textbf{localized by the sequence}  $\pta \Omega_k, X^k\ptc_{k\in
  \mathbb{N}^*},$ if $P\pta \cup_{k=0}^{+\infty}\Omega_k\ptc=1,$
  $\Omega_h \subseteq \Omega_k,$ if $h\leq k,$ and
  $I_{\Omega_k}X^k=I_{\Omega_k}X,$ almost surely for every $k$ in $\mathbb{N}.$
\end{defi}
\begin{rema}  \label{r4.9}
Let $(X_t, 0\leq t\leq T)$ and $(Y, 0\leq t \leq T)$ be two stochastic processes. The following statements are true.
\begin{enumerate}
\item
Let $Y$ and $X$ be localized by the
sequences $\pta \Omega_k, X^k\ptc_{k\in
  \mathbb{N}}$ and $\pta \Omega_k, Y^k\ptc_{k\in
  \mathbb{N}}$, respectively, such that $Y^k$ is $X^k$-forward integrable on $[0,T]$  for every $k$ in $\mathbb{N}.$ Then $Y$ is $X$-forward integrable on $[0,T]$ and 
$$
\int_0^\cdot Yd^-X=\int_0^\cdot
Y^kd^-X^k, \quad \mbox{ on } \Omega_k, \quad a.s..
$$
\item If $Y$ is $X$-forward integrable on $[0,T],$ then $YI_{[0,t]}$ is $X$-forward integrable for every $0\leq t\leq T,$ and 
$$
\int_{0}^\cdot Y_sI_{[0,t]}d^-X_s=\int_0^{\cdot \wedge t} Y_sd^-X_s. 
$$
\item If the covariation process $[X,Y]$ exists on $[0,T],$ then the covariation process $[XI_{[0,t]},YI_{[0,t]}]$ exists for every $0\leq t\leq T,$ and 
$$
\pqa X_{[0,t]},YI_{[0,t]}\pqc=\pqa X,Y\pqc_{t\wedge T}. 
$$      
\end{enumerate}
 \end{rema}

\begin{defi} \label{d2.2}
Let $X=(X_t,0\leq t \leq T)$ and $Y=(Y_t,0\leq t <T)$ be processes with paths respectively in $C^0([0,T])$ and $L^1_{loc}([0,T)),$ i.e. $\int_0^t \vaa Y_s\vac ds< +\infty$ for any $t<T$.
\begin{enumerate}
\item If $YI_{[0,t]}$ is $X$-forward integrable for every $0\leq t < T,$ $Y$ is said \textbf{locally $X$-forward integrable on $[0,T)$}. In this case there exists a continuous process, 
which coincides, on every compact interval $[0,t]$ of $[0,1),$ with the forward integral of $YI_{[0,t]}$ with respect to $X.$ That process will still be denoted with  $I(\cdot,Y,X)=\int_0^\cdot Yd^-X.$
\item If $Y$ is locally $X$-forward integrable and  
$
\lim_{t\rightarrow T} I(t,Y,X)
$ exists almost surely, $Y$ is said $X$-\textbf{improperly forward integrable} on $[0,T]$. 
\item If the covariation process $[X,YI_{[0,t]}]$ exists, for every $0\leq t<T,$  we say that the \textbf{covariation process $[X,Y]$ exists locally} on $[0,T)$ and it is still denoted by $[X,Y].$ In this case there exists a continuous process, 
which coincides, on every compact interval $[0,t]$ of $[0,1),$ with the covariation process $\pqa X,YI_{[0,t]}\pqc.$ That process will still be denoted with  $\pqa X,Y\pqc.$ If $X=Y,$ $\pqa X,X\pqc$ we will say that the \textbf{quadratic variation of $X$ exists locally} on $[0,T].$ 
\item If the covariation process $[X,Y]$ exists locally on $[0,T)$ and $\lim_{t \rightarrow T}[X,Y]_t$ exists, the limit will be called the \textbf{improper covariation process} between $X$ and $Y$  and it will still be denoted by $[X,Y].$
If $X=Y,$ $\pqa X,X\pqc$ we will say that the \textbf{quadratic variation of $X$ exists improperly} on $[0,T].$    
\end{enumerate}
\end{defi}

\begin{rema} \label{r2.5} Let $X=(X_t, 0\leq t\leq T)$ and $ Y=(Y_t, 0\leq t\leq  T)$ be two  stochastic processes being in $C^0([0,1])$ and $L^1([0,1]),$ respectively. 
If $Y$ is $X$-forward integrable on $[0,T]$ then its restriction to $[0,1)$ is $X$-improperly forward integrable and the improper integral coincides with the forward integral of $Y$ with respect to $X.$ 
\end{rema}

\begin{defi} A vector $\pta \pta X^1_t,..., X^m_t\ptc,0\leq t\leq T\ptc $ of continuous processes is said to have all its \textbf{mutual brackets} on $[0,T]$ if $\pqa X^i,X^j\pqc$ exists on $[0,T]$ for every $i,j=1,...,m$.
\end{defi}

In the sequel if $T=1$ we will omit to specify that objects defined above exist on the interval $[0,1]$ (or $[0,1),$ respectively).

\begin{prop}\label{p3.4}
Let  $M=(M_t, 0\leq t\leq T)$ be a continuous
  local martingale with respect to some filtration $\mathbb{F}=\pta \mathcal{F}_t\ptc_{t\in[0,T]}$ of
  $\mathcal{F}.$ Then the following properties hold.
  \begin{enumerate} 
  \item The process $M$ is a finite quadratic variation process on $[0,T]$ and its quadratic variation coincides with the classical bracket appearing in the Doob decomposition of $M^2.$
\item Let $Y=(Y_t,0\leq t\leq T)$ be an $\mathbb{F}$-adapted process with left continuous and bounded paths. Then $Y$ is $M$-forward integrable on $[0,T]$
and $\int_0^\cdot Yd^-M$ coincides with the
classical It\^o integral $\int_0^\cdot YdM.$
\end{enumerate} 
\end{prop}

\begin{prop}\label{p3.6}
Let $V=(V_t,0\leq t\leq T)$  be a  bounded
  variation process and  $
Y=(Y_t,0\leq t\leq T),$ be a process with paths being bounded and with at most countable discontinuities. Then the following properties hold.
\begin{enumerate} 
\item The process $Y$ is $V$-forward
integrable on $[0,T]$ and $\int_0^\cdot Yd^-V$ coincides with the
Lebesgue-Stieltjes integral denoted with $\int_0^\cdot YdV.$
\item The covariation process $\pqa Y,V\pqc$ exists on $[0,T]$ and it is equal to zero. In particular a bounded variation process has zero quadratic variation.
\end{enumerate}
\end{prop}

\begin{coro} \label{c2.15}
Let $X=(X_t,0\leq t\leq T)$ be a continuous process and  $Y=(Y_t,0\leq t \leq T)$ a bounded variation process. Then 
$$
XY-X_0Y_0=\int_0^\cdot X_sdY_s+ \int_0^\cdot Y_sd^-X_s.
$$
\end{coro}

\begin{prop} \label{C1stability} Let $X=(X_t,0\leq t\leq T)$ be a continuous finite quadratic variation process, and $f$ a function in $C^1(\mathbb{R}).$ Then $Y=f(X)$ has a finite quadratic variation on $[0,T]$ and 
$
\pqa Y\pqc=\int_0^\cdot f'(X)^2 d\pqa X\pqc.
$
\end{prop}

\begin{prop} \label{p2.1} Let $X=(X_t,0\leq t\leq T)$ be a continuous finite quadratic variation process and $V=((V^1_t,\dots,V_t^m), 0\leq t\leq T)$ be a vector of continuous bounded variation processes. Then for every $u$ in  $C^{1,2}(\mathbb{R}^m \times \mathbb{R}),$ the process $\pta  \partial_xu(V_t,X_t), 0\leq t\leq T\ptc$ is $X$-forward integrable on $[0,T]$ and
\beqn \nonu
u(V,X)&=&u(V_0,X_0)+\sum_{i=1}^m \int_0^\cdot {\partial_{v_i}} u(V_t,X_t)dV^i_t+\int_0^\cdot \partial_xu(V_t,X_t)d^-X_t\\ \nonu &+&\frac{1}{2}\int_0^\cdot \partial_{xx}^{(2)}u(V_t,X_t)d\pqa X\pqc_t.
\eeqn

\end{prop}
\begin{lemm} \label{l2.1}Let  $X=(X^1_t,\dots,X^m_t, 0\leq t\leq T)$  be a vector of continuous processes having all its mutual brackets. Let  $\psi:\mathbb{R}^m\rightarrow \mathbb{R}$ be of class $C^2(\mathbb{R}^m)$
and $Y=\psi(X).$ 
Then $Z$ is $Y$-forward integrable on $[0,T],$ if and only if $Z\partial_{x^i}\psi(X)$ is $X^i$-forward integrable on $[0,T],$ for every $i=1,...,m$ and
\beqn \nonu 
\int_0^\cdot Zd^-Y=\sum_{i=1}^m \int_0^\cdot Z\partial_{x^i}\psi(X)d^-X^i+\frac{1}{2}\sum_{i,j=0}^m \int_0^\cdot Z\partial_{x^i x^j}^{(2)}\psi(X)d\pqa X^i,X^j\pqc.
\eeqn
\end{lemm}
\begin{proof}
The proof derives from proposition $4.3$ of \cite{RV2}. The result
is a slight modification of that one. It should only be noted that
there forward integral of a process $Y$ with respect to a process $X$ was defined as limit ucp 
of its regularizations.
\end{proof}

\section{Existence of forward integrals and related properties: some examples}\label{s3}
In this section we illustrate examples of processes for which forward integrals exist and we list some related properties which will be extensively used in  further applications to finance. 
\subsection{Forward integrals  of It\^o fields}  \label{s4.0} 
In this subsection $\xi$ will be a $\mathbb{G}$-adapted process with finite quadratic variation, where $\mathbb{G}$ is some filtration of $\mathcal{F}.$
The following definitions and results  are extracted from \cite{FR}.

\begin{defi} \label{Ito-field}
Let $k$ be in $\mathbb{N}^*$. A random field $\pta H(t,x), 0\leq t \leq 1, x \in \mathbb{R}\ptc$ is called a ${C}^{{k}}$ $\mathbb{G}$-{\it{\textbf{It\^{o}-semimartingale field} driven by the vector}} $N=\pta N^1,...,N^n\ptc,$ if $N$ is a vector of semimartingales with respect to $\mathbb{G},$ and   
\beqn \label{IM-field}
H(t,x)=f(x)+\sum_{i=1}^n \int_{0}^t a^i(s,x)dN_s^i, \quad \interval, 
\eeqn
where
$f: \Omega\times \mathbb{R}\rightarrow \mathbb{R}$ belongs to $C^k(\mathbb{R})$ almost surely and it is $\mathbb{G}^0$-measurable for every $x,$
$H$ and $a^i:[0,1]\times \mathbb{R}\times \Omega \rightarrow \mathbb{R},\  i=1,...,n$ are $\mathbb{G}$-adapted  for every $x,$ almost surely continuous with their partial derivatives with respect to $x$ in $(t,x)$ up to order $k,$ and 
for every index $h\leq k$ it holds
$$
\partial_{x}^{(h)} H(t,x)=\partial_x^{(h)}f(x)+\sum_{i=1}^n\int_{0}^t \partial_{x}^{(h)} a^i(s,x)dN_s^i, \quad \interval.
$$
\end{defi}
\begin{defi}\label{d5.14}
We denote with $\mathcal{C}^k_\xi(\mathbb{G})$ the set of processes of the form $$
\pta H(t,\xi_t), 0\leq t \leq 1 \ptc,$$ being $\pta H(t,x), 0\leq t \leq 1, x \in \mathbb{R}\ptc$  a ${C}^{{k}}$ $\mathbb{G}$-{\it{{It\^{o}-semimartingale field} driven by the vector}} $N=\pta N^1,...,N^n\ptc,$ such that $\pta N^1,...,N^n,\xi\ptc$ has all its mutual brackets.  
\end{defi}

\begin{rema} \label{r3.3}
\begin{enumerate}
\item
The set $\mathcal{C}^1_\xi(\mathbb{G})$ is an algebra.
\item Let $\psi$ be in $C^{\infty}(\mathbb{R})$ and $h$ in $\mathcal{C}^2_\xi(\mathbb{G}).$ It\^o formula implies that $\psi(h)$ belongs to  $\mathcal{C}^2_\xi(\mathbb{G}).$ 
\end{enumerate} 
\end{rema}

\begin{prop} \label{p5.14} 
Let $h$ and $k$ be in $\mathcal{C}_\xi^1(\mathbb{G}).$ Then the following statements are true.
\begin{enumerate}
\item  The process $h$ is $\xi$-forward integrable, the forward integral $\int_0^\cdot h_td^-\xi_t$ is the limit ucp of its regularizations and it belongs to $\mathcal{C}_\xi^2(\mathbb{G}).$ 
\item The covariation process  
$
\pta \pqa \int_0^\cdot h_td^-\xi_t, \int_0^\cdot k_td^-\xi_t\pqc \ptc
$ exists and it is equal to 
$   
\int_0^\cdot  h_t k_t d\pqa  \xi\pqc_t.
$  
\item The process $\int_0^\cdot  h_td^-\xi_t$ is forward integrable with respect to the process $\int_0^\cdot k_td^-\xi_t$ and 
$$
\int_0^\cdot h_td^-\int_0^t k_sd^-\xi_s=\int_0^\cdot h_t k_t d^-\xi_t.
$$ 
\end{enumerate}
\end{prop}

Using remark \ref{r4.9} it is not difficult to prove that proposition \ref{p5.14} extends to processes which are \textit{simple combinations} of processes in $\mathcal{C}_\xi^{1}(\mathbb{G}).$ We illustrate this result below.

\begin{defi}\label{d5.10}
Let $\mathcal{S}(\mathcal{C}_\xi^{k}(\mathbb{G}))$ be the set of all processes $h$ of type 
$
h=h^0I_{\{0\}}+\sum_{i=1}^{m}h^iI_{(t_{i-1},t_{i}]}
$
where
$0=t_0\leq t_1,\cdots,t_m=1,$ and $h^i$ belongs to $\mathcal{C}^k_\xi(\mathbb{G}),$ for $i=1,...,m.$
\end{defi}
\begin{rema} \label{r3.8}
Thanks to remark \ref{r3.3}, if $h$ belongs to $\mathcal{S}(\mathcal{C}_\xi^{k}(\mathbb{G}))$ and $\psi$ is of class $C^\infty(\mathbb{R}),$ then $\psi(h)$ is still in $\mathcal{S}(\mathcal{C}_\xi^{k}(\mathbb{G})).$
\end{rema}
\begin{prop}\label{3.8}
Let $h$ and $k$ be in $\mathcal{S}(\mathcal{C}_\xi^{1}(\mathbb{G})).$ Then the we can state the following.  
\begin{enumerate}
\item The process $h$ is $\xi$-forward integrable and it belongs to $\mathcal{S}(\mathcal{C}_\xi^{2}(\mathbb{G})).$
\item The covariation process  
$
\pta \pqa \int_0^\cdot h_td^-\xi_t, \int_0^\cdot k_td^-\xi_t\pqc \ptc
$ exists and it is equal to   
$\int_0^\cdot  h_t k_t d\pqa  \xi\pqc_t.
$  
\item The process $\pta \int_0^\cdot  h_td^-\xi_t, 0\leq t\leq 1 \ptc$ is forward integrable with respect to the process $\pta \int_0^\cdot k_td^-\xi_t, 0\leq t\leq 1 \ptc$ and 
$$
\int_0^\cdot h_td^-\int_0^t k_sd^-\xi_s=\int_0^\cdot h_t k_t d^-\xi_t.
$$ 
\end{enumerate} 
\end{prop}
\begin{proof}
By linearity of forward integral and bilinearity of covariation it is sufficient to prove the statement for processes of type $hI_{[0,t]}$ and $kI_{[0,t]},$ with $h$ and $k$ in $\mathcal{C}^1_\xi(\mathbb{G})$ and $\interval.$ The proof is a consequence of remark \ref{r4.9} and proposition \ref{p5.14}.  
\end{proof}

\subsection{Forward integrals via Malliavin calculus} \label{s4.1}
We work in the Malliavin calculus framework. To this extent we recall some basic notations and definitions from \cite{NP} and \cite{N}.

We suppose that $\pta \Omega, \mathbb{F}, \mathcal{F},P\ptc$ is the canonical probability space, meaning that $\Omega=C\pta [0,1],\mathbb{R}\ptc$, $P$ is the Wiener measure, $W$ is the Wiener process, $\mathbb{F}$ is the filtration generated by $W$ and the $P$-null sets and $\mathcal{F}$ is the completion of the Borel $\sigma$-algebra with respect to $P.$

Let   
$\mathcal{S}$ be the space of all random variables on $(\Omega,\mathcal{F},P),$ of the form
$$
F=f(W(t_1),...,W(t_n)), \quad 0\leq t_0,\cdots,t_n\leq 1,
$$
with $f$ in $C^{\infty}(\mathbb{R}^n)$ being bounded with its derivatives of all orders. The iterated derivative of order $k$ operator is denoted by $D^k.$ Then $D^k:$ $\mathbb{D}^{k,p}\rightarrow L^p\pta\Omega\times[0,1]^k\ptc,$ where $\mathbb{D}^{k,p},$ $p\geq 2,$ $k\in \mathbb{N}^*,$ is the closure of $\mathcal{S}$  with respect to the norm
$$
\vaa \vac F \vaa\vac_{\mathbb{D}^{k,p}}^p=\vaa  \vaa F\vac\vac_{L^p(\Omega)}^p+ \sum_{j=1}^k \vaa \vaa  \vaa \vaa D^jF\vac \vac _{L^2([0,1]^j)}\vac \vac_{L^p(\Omega)}^p.
$$
For any $p\geq 2,$ $L^{1,p}$ denotes the space of all functions $u$ in $L^p\pta \Omega \times [0,1] \ptc$ such that $u_t$ belongs to $\mathbb{D}^{1,p}$ for every  $\interval$  and there exists a measurable version of $\pta D_s u_t, 0\leq s,t \leq 1\ptc$ with
$
\int_0^{1} \mathbb{E}\pqa \vaa \vaa Du_t\vac\vac_{L^{2}([0,1])}^p \pqc dt<\infty.
$
For every $u$ in ${L^{1,p}}$ we denote
$
\vaa \vaa u \vac\vac_{L^{1,p}}^p=\int_0^1 \vaa \vaa u_t \vac\vac_{\mathbb{D}^{1,p}}^pdt.
$
Similarly, for $p\geq 2,$ $L^{2,p}$ denotes the space of all functions
 $u$ in $L^p\pta \Omega\times[0,1] \ptc,$ such that $u_t$ belongs to
 $\mathbb{D}^{2,p}$ for every  $\interval$ and  there exist
 measurable versions of  $\pta D_s u_t, 0\leq s,t \leq 1\ptc$  and $\pta D_r D_s u_t, 0\leq s,t,r\leq 1\ptc$ with
$$
\int_0^{1} \mathbb{E}\pqa \vaa \vaa Du_t\vac\vac_{L^{2}([0,1])}^p \pqc+  \mathbb{E}\pqa \vaa \vaa D^2u_t\vac \vac_{L^{2}([0,1]^2)}^p\pqc dt <\infty.
$$
For every $u$ in $L^{2,p}$ we denote
$
\vaa \vaa u \vac\vac_{L^{2,p}}^p=\int_0^1 \vaa \vaa u_t \vac\vac_{\mathbb{D}^{2,p}}^pdt.
$

The Skorohod integral $\delta$ is the adjoint of the derivative operator $D;$ its domain  is denoted by $Dom\delta.$ An element $u$ belonging to $Dom\delta$ is said Skorohod integrable. We recall that $\mathbb{D}^{1,2}$ is dense in $L^2(\Omega),$ $L^{1,2}\subset Dom\delta,$ and that if $u$ belongs to $L^{1,2}$ then, for each $\interval,$ $uI_{[0,t]}$ is still in  $L^{1,2}.$ In particular it is Skorohod integrable. We will use the notation
$
\delta\pta uI_{[0,t]}\ptc=\int_0^t u_s \delta W_s,
$
for each  $u$ in $L^{1,2}.$ The process $\pta \int_0^t u_s \delta W_s, \interval \ptc$
is mean square continuous and then it admits a continuous version, which will be still denoted by $\int_0^\cdot u_t \delta W_t.$
We finally recall that for every $u$ in $L^{1,p}$ there exists a positive constant $c_p$ such that
\begin{eqnarray}\label{normaskorohod}
\vaa \vaa \delta(u) \vac\vac_{L^p(\Omega)}^p &\leq& c_p \pqa \pta \int_0^1
\vaa \mathbb{E}\pqa u_t  \pqc \vac^2 dt \ptc^{\frac{p}{2}} +
\vaa \vaa \vaa \vaa Du \vac \vac_{L^2([0,1]^2)} \vac\vac_{L^p(\Omega)}^p \pqc \\ \nonu
&\leq & c_p  \vaa \vaa u \vac\vac_{L^{1,p}}^p.
\end{eqnarray}

It is useful to remind the following result contained in \cite{N}, exercise 1.2.13. 
\begin{lemm} \label{4.4} Let $F$ and $G$ be two random variables in $\mathbb{D}^{1,2}.$ Suppose that $G$ and $\vaa \vaa DG\vac \vac_{L^2([0,1])}$ are  bounded. Then $FG$ is still in $\mathbb{D}^{1,2}$ and $D(FG)=FDG+GDF.$  
\end{lemm}

\begin{rema}\label{r4.1}
\begin{enumerate}
\item
Let $u$ be a process in $L^{1,p},$ for some $p\geq 2,$ and $v$  in $L^{1,2}$  such that the random variable
$$
\sup_{t\in[0,1]} \pta \vaa v_t\vac+ \int_0^1  (D_s v_t)^2  ds \ptc
$$
is bounded.
By lemma \ref{4.4} the process $uv$  belongs to $L^{1,p}$ and  $Duv=uDv+vDu.$
\item Let $u$ be a process in $L^{2,p},$ for some $p\geq 2$ and $v$ in $L^{2,2}$ such that the random variable
$$
\sup_{t\in[0,1]} \pta \vaa v_t\vac+ \int_0^1 (D_s v_t)^2 ds+\int_0^1 \int_0^1 (D_r D_s v_t)^2 dr ds \ptc 
$$
is bounded.
Then the process $uv$  belongs to $L^{2,p}.$
\end{enumerate}
\end{rema}
In order to state a chain rule formula we will need the \textit{Fubini}-type lemma below.

\begin{lemm}\label{fubini}
Let $u$ be in $L^2\pta\Omega\times [0,1]^2\ptc.$ Assume that for every $0\leq t\leq 1,$ the process $u(\cdot,t)$ belongs to $L^{1,2},$ that there exist measurable versions of the two processes
$\pta \delta(u(\cdot,t), 0 \leq t\leq 1\ptc$ and $\pta D_r u(s,t), 0\leq r,s,t\leq 1\ptc$  and that 
\beqn \label{i8}
\mathbb{E}\pqa \int_0^1  \vaa \vaa Du(\cdot,t)\vac \vac_{L^{2}([0,1]^2)}^2dt\pqc < +\infty.
\eeqn
Then the process $\pta \int_0^{1}u(s,t)dt, 0\leq s \leq 1)\ptc$ belongs to $L^{1,2}$ and
$$
\delta\pta\int_0^{1}u(\cdot,t)dt\ptc=\int_0^{1}\delta(u(\cdot,t))dt.
$$
\end{lemm}

\begin{proof}
Consider the process  $\pta g_s, 0\leq s\leq 1 \ptc $  so defined:  $g_s=\int_0^1u(s,t)dt.$ Let $0\leq s\leq 1$ be fixed. Since $\pta u(s,t), 0\leq t\leq 1 \ptc$ is in  $L^{1,2},$ $g_s$ is in $\mathbb{D}^{1,2}$ and $Dg_s=\int_0^1 Du(s,t)dt.$ By Fubini theorem $\pta \int_0^1 D_ru(s,t)dt, 0\leq r,s\leq 1\ptc$ admits a measurable version. Thanks to inequality (\ref{i8}),
$\int_0^1 \mathbb{E}\pqa \vaa \vaa Dg_s \vac \vac_{L^2([0,1])}\pqc ds  <+\infty.$ This implies that $g$ is in $L^{1,2}.$ The conclusion of the proof is achieved using exercise 3.2.8, page 174  of \cite{N}. \end{proof}

\begin{defi}\label{d3.10} For every $p\geq 2,$ $L_{-}^{{1,p}}$ will be the space of all processes $u$ belonging to  $L^{1,p}$ such that
$
\lim_{\ep \rightarrow 0}D_tu_{t-\ep}
$
exists in $L^{p}(\Omega \times [0,1])$. The limiting process will be
denoted by $\pta D^-_tu_t, \interval \ptc.$
\end{defi}

\begin{rema}\label{r4.3}
\begin{enumerate}
    \item  If $u$ belongs to $L_{-}^{1,p}$ then
    \beqn \label{f5} \mathbb{E}\pqa \int_0^1 \pta \frac{1}{\ep}\int_{s}^{s+\ep} \vaa D_ru_s-D^-_ru_r \vac^p dr \ptc ds\pqc
    \eeqn
    converges to zero when $\ep$ tends to zero.
Indeed term (\ref{f5}) equals $\frac{1}{\ep}\int_0^\ep f(z)dz,$ with 
$f(z)=\mathbb{E}\pqa \int_0^1 \vaa D_ru_{r-z}-D^-_ru_r \vac^p dr\pqc,$ and $
\lim_{\nu \rightarrow 0}f(z)=0.$
\item Let $u$ and $v$ be two left continuous processes respectively  in $L_{-}^{1,p}$ and $L_{-}^{1,2}$ with $p\geq  2.$  Suppose, furthermore, that $\sup_{t\in[0,1]}\vaa u_t\vac$ belongs to $L^p(\Omega)$ and that the random variable
$
\sup_{t\in [0,1]}\pta \vaa v_t\vac+ \sup_{s\in[0,1]}\vaa D_s v_t\vac \ptc
$
is bounded.
Then $uv$ belongs to $L_{-}^{1,q},$ for every $2 \leq q < p.$ Moreover $D^-uv=uD^-v+vD^-u.$  In particular $v$ belongs to $L_-^{1,q},$ for every $q\geq 2.$
\end{enumerate}
\end{rema}

The hypothesis on the left continuity of $u$ and $v$ on point 2. of previous remark allows us to show that
\beqn \label{f6}
\lim_{\ep \rightarrow 0} \pqa \int_0^1 \vaa z_{t-\ep}-z_t\vac^p dt\pqc=0,  \quad z=u,v, \quad a.s.
\eeqn
That condition could be relaxed. It would be enough to suppose that,  $$\lambda(\pta  0\leq t\leq 1, s.t. \vaa z_t-z_{t^-}\vac\neq 0 \ptc)=0,$$ almost surely, for $z=u,v,$ being $\lambda$ the Lebesgue measure on $\mathcal{B}([0,1])$. 
Nevertheless, convergence in (\ref{f6}) does not hold for every bounded process. To see this it is sufficient to consider, for instance, $z=I_{\mathbb{Q}\cap [0,1]}.$

\begin{lemm} \label{p4.2} Let $u$ and $v$ be respectively in $L_{-}^{1,p},$ $p\geq 2,$ and $L_{-}^{1,2}.$  Suppose that the random variable
$\sup_{t\in [0,1]}\pta \vaa v_t\vac+ \sup_{s\in[0,1]}\vaa D_sv_t\vac \ptc$
is bounded. Then the sequence of processes
$$
\pta \frac{1}{\ep}\int_0^1 u_t v_t (W_{t+\ep}-W_t) dt- \frac{1}{\ep}\int_0^1 u_t \pta \int_t^{t+\ep}v_s\delta W_s \ptc dt \ptc_{\ep>0}
$$
converges in $L^q(\Omega)$ to $\int_0^1 u_t D_t^-v_t dt,$ for every $2 \leq q<p.$
\end{lemm}
\begin{proof}
Set
$$
A_\ep=\frac{1}{\ep}\int_0^1  u_tv_t  \pta W_{t+\ep}-W_t\ptc dt, \quad B_\ep=\frac{1}{\ep}\int_0^1u_t \pta \int_t^{t+\ep}v_s \delta W_s \ptc dt.
$$
Proposition 1.3.4 in section 1.3 of \cite{N} permits to rewrite $A_\ep$ in the following way:
$$
A_\ep=\frac{1}{\ep}\int_0^1  u_tv_t  \pta \int_t^{t+\ep}I_{[0,1]}(s) \delta W_s\ptc dt.
$$
Moreover, by point 1. of remark  \ref{r4.1}, $Duv=vDu+uDv.$ 
For every $0\leq t\leq 1,$ the random variables $\int_t^{t+\ep}D_s(u_t v_t)ds$ and $\int_t^{t+\ep}D_su_t v_sds$ are square integrable. Therefore
property  (4) in section 1.3 of \cite{N} can be exploited to write \beqn \nonu
A_\ep=\frac{1}{\ep}\int_0^1 \pta  \int_t^{t+\ep} u_tv_t  \delta W_s
\ptc dt + \frac{1}{\ep}\int_0^1 \pta \int_t^{t+\ep}D_s(u_tv_t) ds
\ptc dt, \eeqn and \beqn \nonu B_\ep=\frac{1}{\ep}\int_0^1  \pta
\int_t^{t+\ep} u_t v_s  \delta W_s\ptc  dt+ \frac{1}{\ep}\int_0^1
\pta \int_t^{t+\ep} D_s u_t v_sds \ptc dt. \eeqn
This implies
\beqn \nonu
A_\ep-B_\ep&=&\int_0^1 \pta  \frac{1}{\ep}\int_t^{t+\ep} u_t (v_t-v_s)
\delta W_s \ptc dt \\ \nonu &+&\int_0^1  \pta
\frac{1}{\ep}\int_t^{t+\ep}(v_t-v_s)D_su_t ds \ptc dt+
\int_0^1 \pta \frac{1}{\ep} \int_t^{t+\ep}u_t D_sv_t  ds \ptc dt\\ \nonu
&=&I_\ep^1+I_\ep^2+I_\ep^3.
\eeqn
We observe that the function $(\omega,s,t) \longmapsto I_{(t,t+\ep]}(s)
    u_t(\omega)(v_t-v_s)(\omega),$ for every $(\omega,s,t)$ in $\Omega \times [0,1]^2,$ satisfies the hypotheses of lemma
    \ref{fubini}. Therefore $I_\ep^1$ can be rewritten as  follows 
$$
I_\ep^1=\int_0^1 \pta \frac{1}{\ep} \int_{s-\ep}^{s} u_t(v_t-v_s)dt \ptc \delta W_s.
$$
Using inequality (\ref{normaskorohod}) 
it is possible to prove that there exists a positive constant $c$  such that 
$$
\mathbb{E}\pqa \vaa I_\ep^1 \vac ^p\pqc \leq
c \mathbb{E}\pqa \int_0^1  \pta \vaa u_t\vac^p+\pta \int_0^1 (D_r u_t)^2 dr\ptc ^{\frac{p}{2}}\ptc h^\ep_t dt \pqc,
$$
with 
$$
h_t^\ep=\frac{1}{\ep}\int_t^{t+\ep}\pta \vaa v_t-v_s \vac^p+ \pta \int_0^1 \vaa D_rv_t-D_rv_s\vac^2 dr\ptc^{\frac{p}{2}}\ptc ds.  
$$
Since $\sup_{t\in[0,1]}\pta \vaa v_t \vac+\sup_{s\in[0,1]} \vaa D_sv_t\vac\ptc $ is a bounded random variable, for almost all $(\omega,t),$ $h^\ep_t$
converges to zero when $\ep$ goes to zero.
Consequently, Lebesgue dominated convergence theorem applies to conclude that $\mathbb{E}\pqa \vaa I_\ep^1 \vac ^p\pqc$ converges to zero.

Considering the term $I_\ep^2,$ H\"older inequality and the boundedness of $v$ lead to
\beqn \nonu
\mathbb{E}\pqa  \vaa I_\ep^2\vac ^p\pqc 
&\leq& c\mathbb{E}\pqa \int_0^1 \pta \frac{1}{\ep}\int^{t+\ep}_t \vaa D_su_t-D_s^-u_s\vac^p ds \ptc dt \pqc \\ \nonu
&+&c\mathbb{E}\pqa  \int_0^1 \pta \frac{1}{\ep}\int_{s-\ep}^s \vaa v_t-v_s \vac^p dt \ptc  \vaa D_s^-u_s\vac^p ds\pqc,  
\eeqn
for some positive constant $c.$ The first term of previous sum converges to zero by point 1. of remark \ref{r4.3}; the second by Lebesgue dominated convergence theorem.

Finally $I_\ep^3$ may be rewritten as  follows:
\beqn  \nonu
I_\ep^3&=&\int_0^1  \pta \frac{1}{\ep}\int^{t+\ep}_t  (D_s v_t-D_s^- v_s)ds \ptc u_t dt \\ \nonu &+&\int_0^1  \frac{1}{\ep}\int^{t+\ep}_t D_s^-v_s (u_t-u_s)ds dt+ \int_0^1  u_sD_s^-v_sds.
\eeqn
H\"older inequality and again remark \ref{r4.3} implies the convergence to zero in $L^q(\Omega)$ of the first term of the sum for every $2\leq q <p.$
The convergence to zero of the second term of the sum in $L^p(\Omega)$ is due to the boundedness of $\vaa D^-v\vac$ and the following maximal inequality contained in \cite{S}, theorem 1.:
$$
\int_0^1 \sup_{\ep >0}\pta \frac{1}{\ep}\int_{t}^{t+\ep} \vaa z_s\vac^p ds \ptc dt \leq \int_0^1  \vaa z_t\vac^p  dt, \quad z \in L^p(\Omega \times [0,1]). 
$$
 This leads to the conclusion.
\end{proof}
We  omit the proof of the following lemma which is, indeed, a slight modification of the proof of previous one.
  
\begin{lemm} 
\label{p4.3} Let $v$ be in $L_{-}^{1,p},$ $p\geq 2$. Then the sequence of processes
$$
\pta \frac{1}{\ep}\int_0^1 v_t (W_{t+\ep}-W_t) dt- \frac{1}{\ep}\int_0^1 \pta \int_t^{t+\ep}v_s\delta W_s \ptc dt \ptc_{\ep>0}
$$
converges in $L^p(\Omega)$ to $\int_0^1 D_t^-v_t dt.$
\end{lemm}
\begin{lemm} \label{l4.4}
Let $u$ be a process in $L^{1,p}$ with $p\geq 2.$
Then the process $$\pta \int_0^t u_s ds, \interval\ptc$$ belongs to
$L^{1,p}_{-},$  and $D^-\pta \int_0^\cdot  u_t dt\ptc=\int_0^\cdot Du_t dt.$
\end{lemm}
\begin{proof}
We set $g_t=\int_0^t u_sds.$ Clearly $g$ is in $L^{p}(\Omega \times [0,1]).$ As already observed for the proof of lemma \ref{fubini}, since the process $u$ belongs to $L^{1,2}$ for every $\interval,$ $g_t$ is in $\mathbb{D}^{1,2}$ and $Dg_t=\int_0^t Du_s ds.$ Moreover H\"older inequality implies
\beqn \nonu
\mathbb{E} \pqa \int_0^1 \vaa \vaa Dg_t \vac \vac_{L^2([0,1])}^{p} dt\pqc&\leq&\mathbb{E}\pqa \int_0^1  \vaa \vaa Du_s \vac \vac_{L^2([0,1])}^{p} ds\pqc  <+\infty.
\eeqn
Then $g$ belongs to $L^{1,p}.$
To conclude it is sufficient to observe that
$$
\mathbb{E} \pqa  \int_0^1  \int_{t-\ep}^t \vaa D_t u_s\vac^p ds  dt\pqc=\mathbb{E} \pqa  \int_0^1  \int^{s+\ep}_s \vaa D_t u_s\vac^p dt  ds\pqc,
$$
and that the right hand side of previous equality  converges to zero when  $\ep$ goes to zero by point $1.$ of remark \ref{r4.3}. 
\end{proof}

\begin{lemm} \label{l4.6}
Let $u$ be a process in $L^{2,p},$ with $p\geq 2.$ Suppose furthermore that 
\beqn \label{f7}
\int_0^1 \pta \mathbb{E}\pqa \vaa \vaa Du_t \vac\vac_{L^p( [0,1])}^p\pqc + \mathbb{E}\pqa \vaa \vaa D^2u_t  \vac\vac_{L^p([0,1]^2)}^p\pqc \ptc dt <+\infty.
\eeqn
Then the process $\pta \int_0^t u_s \delta W_s, \interval \ptc$ is in $L^{1,p}_{-},$ and  $$D^{-}\pta \int_0^\cdot u_t \delta W_t\ptc= \int_0^\cdot D u_t \delta W_t.$$
\end{lemm}

\begin{proof}
We set $g=\int_0^\cdot u_t\delta W_t.$
By proposition 5.5 of \cite{NP}, for every $t$ in $[0,1],$ $g_t$ belongs to $\mathbb{D}^{1,2}$ and
$
D_r g_t=\delta \pta  D_r uI_{[0,t]}\ptc +u_rI_{[0,t]}(r) , 
$ 
for every $r,$ almost surely. 
Using inequality (\ref{normaskorohod}) it is possible to find a positive constant $c$ such that
$$
\vaa \vaa g \vac \vac_{L^p(\Omega \times [0,1])}^p \leq c\vaa \vaa u \vac\vac_{L^{1,p}}^p<+\infty.
$$
To prove that  $g$ belongs to $L^{1,p}$ we still have to show that
$
\mathbb{E}\pqa \int_0^1 \vaa \vaa Dg_t \vac \vac_{L^2([0,1])}^{p} dt\pqc
$ is finite.
Clearly, $\mathbb{E}\pqa \int_0^1 \vaa \vaa u I_{[0,t]}\vac \vac_{L^2([0,1])}^{p} dt\pqc\leq \vaa \vaa u \vac \vac_{L^p(\Omega \times [0,1])},$ which is finite.
It remains to prove that
$
\mathbb{E}\pqa \int_0^1 \pta \int_0^1 \vaa \delta \pta D_r u I_{[0,t]}\ptc \vac^2 dr \ptc^{\frac{p}{2}} dt\pqc <+\infty.
$
Applying again inequality (\ref{normaskorohod}) we obtain, for some $c>0,$
\beqn \nonu
\mathbb{E}\pqa \int_0^1 \pta \int_0^1 \vaa \delta \pta D_r u I_{[0,t]}\ptc \vac^2 dr \ptc^{\frac{p}{2}} dt\pqc &\leq& c \int_0^1  \mathbb{E}\pqa\int_0^1 \vaa \delta \pta D_r u I_{[0,t]}\ptc \vac^p dt \pqc dr\\ \nonu
&\leq & c \int_0^1  \int_0^1 \vaa \vaa  D_ruI_{[0,t]} \vac \vac_{L^{1,p}}^p dr dt  \\ \nonu
&\leq & c \int_0^1 \vaa \vaa  D_ru \vac \vac_{L^{1,p}}^p dr.
\eeqn
Last term in the expression above is bounded by the integral appearing in inequality  (\ref{f7}). This permits to get the result.
\end{proof}

\begin{prop} \label{c4.1} Let $v$ be a process in $L_{-}^{1,p},$ with $p>4$.  Then $v$ is both forward and Skorohod integrable with respect to $W$ and
\beqn \nonu
\int_0^\cdot v_t d^- W_t=\int_0^\cdot v_t \delta W_t+\int_0^\cdot D_t^-v_tdt.
\eeqn
Furthermore, if $v$ is also left continuous with right limit, then  $\int_0^\cdot v_td^-W_t$ has finite quadratic variation equal to $\int_0^\cdot v_t^2dt.$ 
\end{prop}
\begin{proof}
First of all we observe that if a process $v$ belongs to $L^{1,p}_-$ then $I_{[0,t]}v$ inherits the property for every $t$ in $[0,1].$ 
Using lemma \ref{p4.3} and lemma \ref{fubini} we find that $I(\ep,v,W,t)-\int_0^t v_s\delta W_s$ converges in $L^p(\Omega)$ toward $\int_0^t D_s^-v_s ds,$ for every $0\leq t\leq 1.$ 
If $v$ belongs to $L_{-}^{1,p},$ with $p>4,$ by theorem 5.2 in  \cite{NP}, the Skorohod integral process $\int_0^\cdot v_t\delta W_t$ admits a continuous version. At the same time,  thanks to theorem 1.1 of \cite{RV1} we know that $\int_0^\cdot v_t\delta W_t$ has finite quadratic variation equal to $\int_0^\cdot v_t^2 dt.$ The proof is complete.
\end{proof}

\begin{prop} \label{p4.2a} Let $u$ and $v$ be left continuous processes, respectively in $L_{-}^{1,p}$ and $L_{-}^{1,2},$ with $p>4$. Suppose that  $\sup_{t\in[0,1]}\vaa u_t\vac$ belongs to $L^p(\Omega),$ and that the random variable  
$\sup_{t\in[0,1]} \pta \vaa v_t\vac+ \sup_{s\in[0,1]}\vaa D_sv_t\vac \ptc 
$ is bounded.
Then
$uv$ and $v$ are forward integrable with respect to $W.$ Furthermore $u$ is forward integrable with respect to $\int_0^\cdot v_td^-W_t$ and
\beqn \nonu
\int_0^\cdot u_t d^- \pta \int_0^t v_s d^- W_s \ptc &=&\int_0^\cdot u_t v_t d^-W_t \\ \nonu &=& \int_0^\cdot u_tv_t \delta W_t +\int_0^\cdot (v_tD_t^-u_t+u_tD_t^-v_t) dt.
\eeqn
\end{prop}
\begin{proof}
By point 2. of remark \ref{r4.3} the process $uv$ belongs to $L^{1,q}_-,$ for every $4<q< p,$ and $D^-uv=vD^-u+uD^-v.$ Proposition \ref{c4.1} immediately implies that
$$
\int_0^\cdot u_t v_t d^-W_t=\int_0^\cdot u_tv_t \delta W_t +\int_0^\cdot (v_tD_t^-u_t+u_tD_t^-v_t) dt.
$$
Lemma \ref{p4.2} permits to write, for every $\interval,$
\beqn\nonu
I\pta \ep,u,\int_0^\cdot v_td^-W_t,t\ptc&=&\frac{1}{\ep} \int_0^t u_s\int_s^{s+\ep}v_r d^-W_r ds\\ \nonu &=&\frac{1}{\ep}\int_0^t u_s \pta \int_s^{s+\ep}v_r \delta W_r\ptc ds\\ \nonu
&+&\frac{1}{\ep}\int_0^t u_s \pta \int_s^{s+\ep} D_r^-v_r dr\ptc  ds.
\eeqn
Since $\sup_{t\in[0,1]}{\vaa D^-_tv_t\vac}$ belongs to $L^p(\Omega)$ the second term of previous sum converges toward $\int_0^t u_sD^-_sv_sds$ in $L^q(\Omega),$ for every $2\leq  q<p.$
As a consequence of this, by lemma \ref{p4.2}, $I\pta \ep,u,\int_0^\cdot v_td^-W_t,t\ptc$ converges toward $\int_0^t u_sv_s d^-W_s$ in $L^2(\Omega).$ The proof is then complete.
\end{proof}

\begin{defi}
We say that a process $u$  belongs to $L^{1,p}_{-,loc}$ if it
is localized by a sequence $\pta \Omega_k,
u^k \ptc_{k\in \mathbb{N}},$ with   $u^k$ belonging to
$L^{1,p}_{-} $ for every $k$ in $\mathbb{N}.$
\end{defi}

\begin{lemm} \label{l4.8}
Let $u=\pta u^1,\dots, u^n \ptc,$ $n>1,$ be a vector of left continuous processes with bounded paths and in $L^{1,p}_{-},$ for some $p \geq 2$.
Then, for every $\psi$ in $C^1\pta \mathbb{R}^n\ptc$ the process $\psi(u)$ belongs to $L^{1,p}_{-,loc}$. Moreover, the localizing sequence $\pta \Omega_k, \psi(u)^k\ptc_{k\in
  \mathbb{N}}$ can be chosen such that $\psi(u)^k$ is left continuous,  and $\sup_{t\in[0,1]}\vaa \psi(u)^k_t\vac$ belongs to $L^p(\Omega)$ for every $k$ in $\mathbb{N}.$
\end{lemm}
\begin{proof}
For $k$ in $\mathbb{N}^*,$  set
$
\Omega_k=\pga \sup_{0\leq t\leq 1} ||u_t||_{\mathbb{R}^n} \leq k\pgc
$
and  $\psi(u)^k=\psi(u)f_k(u),
$
being $f_k(u)=f(\frac{u}{k}),$ and $f$ a smooth function from
$\mathbb{R}^n$ to $\mathbb{R},$ with compact support and $f(x)=1,$ for every
$\vaa \vaa x  \vac\vac\leq 1.$
Clearly $\psi(u)$ is localized by  $\pta \Omega_k, \psi(u)^k\ptc_{k\in
  \mathbb{N}}.$ By \cite{NP}, proposition  4.8, $ \psi(u)^k$ belongs to
$L^{1,2},$ for every $k$ in $\mathbb{N}^*.$  Since $\psi \circ
f_k$ has bounded first partial derivatives, proposition 1.2.2 of \cite{N} implies that
$$
D \psi(u)^k_s=\sum_{i=1}^n \partial_i(
\psi\circ f_k) (u_s)D u^i_s.
$$
In particular $\psi(u)^k$ belongs to $L^{1,p}.$
Using the  continuity of all first partial derivatives of  $\psi \circ
f_k$ and the left continuity of $u^i$ for every $i=1,...,n,$ it is possible to prove that $\psi(u)^k$ belongs indeed to  $L^{1,p}_{-},$ and  $D^-\psi(u)^k=\sum_{i=1}^n \partial_i
(\psi\circ f_k (u))D^- u^i.$  The proof is then complete.
\end{proof}

We conclude this section giving a generalization of proposition $\ref{p4.2a}.$

\begin{prop} \label{p4.15} Let $u=\pta u^1,\dots, u^n \ptc,$ $n>1,$ be a vector of left continuous processes with bounded paths and in $L^{1,p}_{-},$ with $p>4.$ Let $v$ be a process in $L_{-}^{1,2}$ with left continuous paths such that the random variable
$
\vaa v_t\vac+\sup_{s\in[0,1]} \vaa D_sv_t\vac
$ 
is bounded.
Then for every $\psi$ in $C^1(\mathbb{R}^n)$
$\psi(u)v$ and $v$ are forward integrable with respect to $W.$ Furthermore $\psi(u)$ is forward integrable with respect to $\int_0^\cdot v_td^-W_t$ and
\beqn \nonu
\int_0^\cdot \psi(u_t)d^- \pta \int_0^t v_s d^- W_s \ptc &=&\int_0^\cdot \psi(u_t) v_t d^-W_t.
\eeqn
\end{prop}
\begin{proof}
Let $\pta \Omega^k, \psi(u)^k\ptc_{k\in \mathbb{N}}$  be
a localizing sequence for $\psi(u)$ such that $\psi(u)^k$ is left continuous and $\sup_{t\in[0,1]}\vaa \psi(u)^k_t\vac$ belongs to $L^p(\Omega)$ for every $k$ in $\mathbb{N}.$ Such a sequence exists thanks to lemma \ref{l4.8}. Clearly $\pta \Omega^k, \psi(u)^kv\ptc_{k\in \mathbb{N}}$ localizes $\psi(u)v.$
For every $k$ in $\mathbb{N},$  thanks to proposition \ref{p4.2a}, $\psi(u)^k$ and $\psi(u)^kv$ are forward integrable with respect to $W$ and
$$
\int_0^\cdot \psi(u)^k_td^-\int_0^t v_sd^-W_s=\int_0^\cdot \psi(u)^k_t v_td^-W_t.
$$
The conclusion follows by remark \ref{r4.9}.
\end{proof}

\subsection{Forward integrals of anticipating processes: substitution formulae}\label{s4.2}
Let $\mathbb{F}=\pta \mathcal{F}_t\ptc_{t\in[0,1]}$ be a filtration on $\pta \Omega, \mathcal{F}, P\ptc,$ with $\mathcal{F}_1=\mathcal{F},$ and $L$ an $\mathcal{F}$ measurable random variable with values in $\mathbb{R}^d$.
We set $
\mathcal{G}_t=\pta \mathcal{F}_{t} \vee \sigma(L)\ptc,
$ and we suppose that $\mathbb{G}$ is right continuous:
$$
\mathcal{G}_t=\bigcap_{\ep > 0}\pta \mathcal{F}_{t+\ep} \vee \sigma(L)\ptc.
$$

In this section $\mathcal{P}^{\mathbb{F}}$ ($\mathcal{P}^{\mathbb{G}},$ resp.) will denote the $\sigma$-algebra of $\mathbb{F}$(of $\mathbb{G},$ resp.)-predictable processes.
$E$ will be the Banach space of all continuous functions on $[0,1]$ equipped with the uniform
norm $\left|\left|f\right|\right|_{E}=\sup_{t\in[0,1]}\left|f(t)\right|.$

\subsubsection{Preliminary results}
We state in the sequel some results about forward integrals involving processes that are random field evaluated at $L.$ To be more precise we will establish conditions to insure the existence of such integrals, their quadratic variation and a related associativity property.
\begin{defi}An increasing sequence of random times $\pta
T_k\ptc_{k\in\mathbb{N}}$ is said \textit{\textbf{suitable}} if $P\pta
\cup_{k=0}^{+\infty}\pga
T_k=1 \pgc \ptc =1.$
\end{defi}

\begin{defi}\label{d4.16} For every random time $0 \leq S \leq 1,$ $p>0,$ and $\gamma>0,$ we define $\texttt{{C}}^{{p},{\gamma}}_{S}$ as the set of all families of continuous processes \textbf{$((F(t,x),\interval);x\in \mathbb{R}^d)$} such that for each compact set $C$ of $\mathbb{R}^d$ there exists a constant $c>0$  such that
$$
\mathbb{E}\pqa \sup_{t\in[0,S]}\left|F(t,x)-F(t,y)\right|^{p}\pqc\leq c \left|x-y\right|^{\gamma}, \quad \forall x, y\in C.
$$
If $S=1,$ $\texttt{C}^{p,\gamma}$ will stand for $\texttt{C}_S^{p,\gamma}$.
\end{defi}

We begin recalling a result stated in \cite{RV1}, lemma 1.2, page 93.
\begin{lemm}
\label{ls}Let $\pga (F_n(t,x), \interval),(F(t,x),\interval); n\geq 1, x\in \mathbb{R}^d\pgc $ be a family  of continuous processes such that $F_n$ and $F$ are $\mathcal{F}\otimes\mathcal{B}([0,1])\otimes\mathcal{B}(\mathbb{R}^d)$-measurable. Suppose that for each $x$ in $\mathbb{R}^d,$ $F_n(\cdot,x)$ converges to $F(\cdot,x)$ ucp and that there exist $p>1,$ $\gamma>d,$ with
$$
\mathbb{E}\pqa \sup_{t\in[0,1]}\left|F_n(t,x)-F_n(t,y)\right|^{p}\pqc\leq c \left|x-y\right|^{\gamma}, \quad \forall x, y\in C, \quad \forall n\in \mathbb{N}.
$$
Then $x\mapsto F(\cdot,x)$ admits a continuous version $\bar{F}(\cdot,x)$, from $\mathbb{R}^d$ to $E$ 
and $F_n(\cdot,L)$ converges toward $F(\cdot,L)$ ucp.
\end{lemm}

\begin{defi}\label{d4.17} ${\texttt{bL}}(\mathcal{{P}}^{\mathbb{{F}}}{\otimes} \mathcal{{{B}}}(\mathbb{{R}}^{d}))$ will denote the set of all functions  $$((h(t,x),\interval);x\in\mathbb{R}^d)$$ which are   $\mathcal{P}^{\mathbb{F}}\otimes\mathcal{B}(\mathbb{R}^d)$-measurable, such that for every $x$ in $\mathbb{R}^d,$ $h(\cdot,x)$ has left continuous and bounded paths.
\end{defi}
\begin{defi}\label{d5.13a}
Let $p>1, \gamma >0$. We define $\mathcal{A}^{p,\gamma}$ as the set of all functions $h$ in  ${\texttt{bL}}(\mathcal{{P}}^{\mathbb{{F}}}{\otimes} \mathcal{{{B}}}(\mathbb{{R}}^{d}))$ satisfying the following assumption. There exists a suitable sequence of stopping times $\pta T_k\ptc_{k\in \mathbb{N}}$ such that $h$ belongs to $\bigcap_{k\in{\mathbb{N}}}\texttt{{C}}^{{p},{\gamma}}_{T_k}.$
\end{defi}

We state this lemma which will be useful later.
\begin{lemm}\label{l4.17}
Let $h$ and $g$ be respectively in $\mathcal{A}^{p,\gamma_p}$ and $\mathcal{A}^{q,\gamma_q}$ for some $p,q>1,$ $\gamma_p,\gamma_q>0.$ Then the following statements hold.
\begin{enumerate}
\item If $\psi$ belongs to $C^1(\mathbb{R})$ and it has bounded derivative, then $\psi(h)$ belongs to $\mathcal{A}^{p,\gamma}.$
\item 
The process $hg$  belongs to  $\mathcal{A}^{\alpha,\frac{\gamma \alpha}{p}\wedge \frac{\gamma_q \alpha}{q}}$ with $\alpha=\frac{pq}{p+q}.$
\item If $N$ is a continuous $\mathbb{F}$-semimartingale, then the function $$\pta \int_0^t h(s,x)dN_s, \interval, x\in \mathbb{R})\ptc$$ belongs to $\mathcal{A}^{p,\gamma}$.
\end{enumerate}
\end{lemm}
\begin{proof}
Let $\pta T_k\ptc_{k\in \mathbb{N}}$ be a suitable sequence of stopping times such that $h$ belongs to $\bigcap_{k\in{\mathbb{N}}}\texttt{{C}}^{{p},{\gamma}}_{T_k}.$ 

The conclusion of the first point is straightforward. 

Concerning the second point we set, for every $k$ in $\mathbb{N,}$
$$
S_k=\inf \pga \interval, |  \vaa h(t,0)\vac+\vaa g(t,0)\vac\geq k\pgc \wedge T_k.
$$
If $C$ is a compact set of $\mathbb{R}^d,$ using H\"older inequality we obtain
\beqn \nonu
\mathbb{E}\pqa \sup_{t\in[0,S_k]}\vaa hg(t,x)-hg(t,y)\vac^{\alpha}\pqc &\leq& c \pta \mathbb{E}\pqa \sup_{t\in[0,S_k]}\vaa h(t,x)-h(t,y)\vac^{{p}}\pqc\ptc^{\frac{\alpha}{p}} \\ \nonu
&+& c \pta \mathbb{E}\pqa \sup_{t\in[0,S_k]}\vaa g(t,x)-g(t,y)\vac^{{q}}\pqc\ptc^{\frac{\alpha}{q}}  \\ \nonu
&\leq & c \vaa x-y\vac^{\frac{\gamma \alpha}{p}\wedge \frac{\gamma_q \alpha}{q}},
\eeqn
where
$c=\sup_{x,y \in C}\pta \mathbb{E}   \pqa  \sup_{t\in[0,S_k]} \vaa h(t,x)\vac^{{p}}   \pqc  \ptc^{\frac{\alpha}{p}}+\pta\mathbb{E}\pqa  \sup_{t\in[0,S_k]} \vaa g(t,x) \vac^{{q}}\pqc\ptc^{\frac{\alpha}{q}}$
is bounded thanks to the choice of the sequence $(S_k)_{k\in \mathbb{N}}$ and the compactness of $C.$

To prove point 3. it is sufficient to define
$
S_k=\inf \pga \interval, | \vaa V\vac_t+ [M]_t \geq k\pgc \wedge T_k
$, for every $k$ in $\mathbb{N},$ where $N=M+V,$ $M$ is an $\mathbb{F}$-local martingale, $V$ is a bounded variation process, and $\vaa V\vac$ denotes the total variation of $V.$  
\end{proof}
\subsubsection{Existence}
Let $M,V:\pta \Omega\times
[0,1]\times\mathbb{R}^d,\mathcal{F}\otimes\mathcal{B}([0,1])\otimes\mathcal{B}(\mathbb{R}^d)\ptc\longrightarrow
\pta \mathbb{R},\mathcal{B}(\mathbb{R})\ptc$ be measurable functions such
that for each $x$ in $\mathbb{R}^d,$ $M(0,x)=V(0,x)=0,$ $(M(t,x),\interval)$ is an
\filtrationf-continuous local martingale, and $(V(t,x), \interval),$ a continuous  bounded variation process.

\begin{rema}\label{r5.2} If $h$ is in $\texttt{bL}(\mathcal{P}^{\mathbb{F}}\otimes \mathcal{B}(\mathbb{R})),$ the process $(h(t,L), \interval)$ is left continuous  with bounded paths and the process $(V(t,L), \interval)$ is continuous with bounded variation. Then, by proposition \ref{p3.6}, $\int_0^\cdot h(t,L)d^-V(t,L)$ exists and coincides with the Lebesgue-Stieltjes integral $\int_0^\cdot h(t,L)dV(t,L).$ Moreover
\beqn \nonu
\int_0^\cdot h(t,L)dV(s,L)=\pta \int_0^\cdot h(t,x)dV(t,x)\ptc_{x=L}.
\eeqn
\end{rema}

\begin{lemm}\label{l5.4} Let $h$ be in $\texttt{bL}(\mathcal{P}^{\mathbb{F}}\otimes \mathcal{B}(\mathbb{R})).$ Suppose that both $\sup_{t\in[0,1]} \vaa h(t,0)\vac $ and $\sup_{t\in[0,1]}\vaa M(t,0)\vac$ are bounded, and that there exist $p>1,$ $q>\frac{p}{p-1}$, $\gamma_p>\frac{d(q+p)}{q},$ $
\gamma_q>\frac{d(q+p)}{p}$ such that $M$ belongs to $\texttt{C}^{p,\gamma_p}$ and $h$ to $\texttt{C}^{q,\gamma_q}.$
Then the function $x \mapsto  \int_0^\cdot h(s,x)dM(s,x)$ admits a continuous version, $\int_0^\cdot h(s,L)d^-M(s,L)$ exists as limit ucp of its regularizations and
\beqn \nonu
 \int_0^\cdot h(t,L)d^-M(t,L)=\pta\int_0^\cdot h(t,x)dM(t,x)\ptc_{x=L}.
\eeqn
\end{lemm}

\begin{proof}
For every $x\in \mathbb{R}^d$ and $\interval$ we set
\beqn \nonu
F_\ep(t,x)=\frac{1}{\ep}\int_0^t h(s,x)(M(s+\ep,x)-M(s,x))ds,
\eeqn
and 
$$
F(t,x)=\int_0^th(s,x)dM(s,x).
$$
To prove our statement we verify that lemma \ref{ls} applies
to the families defined above.

Let $x$ and $y$ in $\mathbb{R}^d$ be fixed. Point 2. of proposition \ref{p3.4} implies
that $F_{\ep}(\cdot,x)$ converges $ucp$ to $F(\cdot,x)$. Set $\alpha=\frac{pq}{p+q}.$ Using theorem  45 in chapter IV of \cite{P},  we can write, for every $\ep>0,$
\beqn \nonu
F_\ep(\cdot,x)-F_\ep(\cdot,y)&=&\int_0^\cdot \pta \frac{1}{\ep}\int_{r-\ep}^r (h(s,x)-h(s,y))ds \ptc dM(r,x)\\ \nonu
&+&\int_0^\cdot \pta \frac{1}{\ep}\int_{r-\ep}^r h(s,y)ds \ptc d\pta M(r,y)-M(r,x)\ptc.\nonu
\eeqn
Thanks to theorem 2 in chapter V of \cite{P}, we find a positive constant $a,$ depending only on $p$ and $q$ such that, for every $\ep>0,$
$$
\mathbb{E}\pqa  \sup_{t\in[0,1]}\vaa
F(\ep,t,x)-F(\ep,t,y)\vac^\alpha \pqc \leq a(
\delta^{}_1+\delta^{}_2)
$$
with
$$
\delta^{}_1=\mathbb{E}\pqa
\sup_{t\in[0,1]} \vaa
h(t,x)-h(t,y)\vac^q\pqc^\frac{\alpha}{q}\mathbb{E}\pqa\sup_{t\in[0,1]} \vaa
M(t,x)\vac^p\pqc^\frac{\alpha}{p}
$$
and
$$
\delta^{}_2=\mathbb{E}\pqa
\sup_{t\in[0,1]} \vaa
M(t,x)-M(t,y)\vac^p\pqc^\frac{\alpha}{p}\mathbb{E}\pqa\sup_{t\in[0,1]} \vaa
h(t,y)\vac^q\pqc^\frac{\alpha}{q}.
$$
Thanks to the hypotheses on $M$ and $h$ 
it is possible to find a constant ${b}$
depending on $C$ such that
\beqn \nonu
\delta^{}_1\leq b \vaa x-y\vac^{\frac{\gamma_q \alpha}{q}},
\quad \delta^{}_2\leq b \vaa x-y\vac^{\frac{\gamma_p \alpha}{p}}, \quad \forall \ep>0.
\eeqn
Consequently, there will exist $c>0$ such that
$$
\mathbb{E}\pqa  \sup_{t\in[0,1]}\vaa
F(\ep,t,x)-F(\ep,t,y)\vac^\alpha \pqc \leq c\vaa
x-y\vac^{\gamma}, \quad \quad \forall x,y \in C, \quad \forall \ep>0,
$$
with $\gamma=\frac{\gamma_p \alpha}{p} \wedge \frac{\gamma_q \alpha}{q}>d$ and
proof is complete.
\end{proof}

The following  proposition represents a generalization of previous lemma.

\begin{prop}\label{p5.5}
Suppose that $M$ belongs to $\mathcal{A}^{p,\gamma_p}$ and that $h$ belongs to $\mathcal{A}^{q,\gamma_q}$ for some $p>1,$ $q>\frac{p}{p-1}$, $\gamma_p>\frac{d(q+p)}{q}$ and $
\gamma_q>\frac{d(q+p)}{p}.$ 
Then $x \mapsto  \int_0^\cdot h(s,x)dM(s,x)$ admits a continuous version, $\int_0^\cdot h(s,L)d^-M(s,L)$ exists as limit ucp of its regularizations and
\beqn \nonu
 \int_0^\cdot h(s,L)d^-M(s,L)=\pta\int_0^\cdot h(s,x)dM(s,x)\ptc_{x=L}.
\eeqn
\end{prop}

\begin{proof}
We observe that we do not loose generality assuming that there exists a suitable sequence of $\mathbb{F}$-stopping times $\pta T_k\ptc_{k\in \mathbb{R}}$ such that $M$ and $h$ belong respectively to $\bigcap_{k\in \mathbb{N}} \texttt{C}^{p,\gamma_p}_{T_k}$ and $\bigcap_{k\in \mathbb{N}}\texttt{C}^{q,\gamma_q}_{T_k}.$ Let $(S_k)_{k\in \mathbb{N}}$ be a suitable sequence of $\mathbb{F}$-stopping times such that, for every $k$ in $\mathbb{M},$ $S^k$ is the first instant, between $0$ and $1,$ the process $\vaa M(\cdot,0)\vac+\vaa h(\cdot,0)\vac$ is greater than $k.$ Set, for every $k$ in $ \mathbb{N},$
$$
R_k=S_k\wedge T_k,\quad M^k=M^{R_k}, \quad h^k=h^{R_k},
$$
and
$$
\Omega_k=\left\{\sup_{ t\in[0,1]}\left|M(t,0)\right|\leq
k \pgc \cap \pga \sup_{t\in[0,1]}\left|h(t,0)\right|\leq k \pgc\cap \pga  R_k=1\pgc.
$$
Let $k$ be fixed. It is clear that $\sup_{t\in[0,1]}\vaa h^k(t,0)\vac$ and $\sup_{t\in[0,1]}\vaa M^k(\cdot,0)\vac$ are bounded
and that $M^k$ and $h^k$ belong, respectively, to $\texttt{C}^{p,\gamma_p}$ and $\texttt{C}^{q,\gamma_q}.$ 
We can thus apply lemma \ref{l5.4} to state that the function $x\mapsto \int_0^\cdot h^k(t,x)dM^k(t,x)$ admits a continuous version and
$$
\int_0^\cdot h^k(t,L)dM^k(t,L)=\pta \int_0^\cdot h^k(t,x)dM^k(t,x)\ptc_{x=L}.
$$
By the \textit{local character} of the classical stochastic integral, see \cite{P}, theorem 26, $F^k(\cdot,x)=F^h(\cdot,x)=\int_0^\cdot h(t,x)dM(t,x),$ for every $x$ in $\mathbb{R}^d,$ almost surely on $\Omega_h,$ for every $h\leq k.$ Therefore it is possible to define $\pta \bar{F}(t,x),\interval, x \in \mathbb{R}^d\ptc$ such that $x \mapsto \bar{F}(\cdot,x)$ is continuous, and for every $k$ in $\mathbb{N}$
$$
\bar{F}(\cdot,x)=\bar{F}^k(\cdot,x)=\int_0^\cdot h(t,x)dM(t,x),\quad \forall x \in \mathbb{R}^d, \quad \mbox{ on } \Omega_k.
$$
Furthermore, remark \ref{r4.9} implies that
$$
\lim_{\ep\rightarrow 0}\frac{1}{\ep}\int_0^\cdot
h(s,L)(M(s+\ep,L)-M(s,x))ds=\bar{F}(\cdot,L),
$$
ucp, since the convergence holds on every on $\Omega_k,$ for every $k$ in  $\mathbb{N}.$
\end{proof}

Previous proposition combined with lemma \ref{l4.17} implies directly the following corollary.

\begin{coro}\label{c5.6}
Let $N$ be a continuous $\mathbb{F}$-local martingale. Let $h$ be in $\mathcal{A}^{q,\gamma_q},$
for $q>1$, $\gamma_q>d.$
Then $x \mapsto  \int_0^\cdot h(t,x)d N_t$ admits a continuous version, $\int_0^\cdot h(t,L)d^-N_t$ exists as limit ucp of its regularizations and
\beqn \nonu
 \int_0^\cdot h(t,L)d^-N_t=\pta\int_0^\cdot h(t,x)dN_t\ptc_{x=L}.
\eeqn
\end{coro}

\subsubsection{Quadratic variation}
We examine the existence of the quadratic variation of forward integrals of the anticipating processes considered in the present subsection. We start giving a generalization of a substitution formula proved in \cite{RV1}, proposition 1.3. We furnish, in fact, a \textit{localized} version of that result in view of further applications to finance. 
We omit the details of its proof, which are indeed similar to those used in the proof of proposition \ref{p5.5}.

\begin{prop}\label{p5.9}
Suppose that $M$ belongs to $\mathcal{A}^{p,\gamma}$ with $p>2$ and $\gamma>2d.$  
Then $x\mapsto \pqa M(\cdot,x), M(\cdot,x) \pqc $ admits a continuous
version, the process $M(\cdot,L)$ has finite quadratic variation and
$$
\pqa M(\cdot,L), M(\cdot,L)\pqc=\pqa M(\cdot,x),M(\cdot,x)\pqc_{x=L}.
$$
\end{prop}
A consequence of previous proposition is the following.

\begin{prop}\label{p5.5  a}
Suppose that $M(t,x)=\int_0^\cdot h(t,x)dN(t,x),$ for every $x$ in $\mathbb{R}^d,$ where $h$ and $N$ verify  the following assumption. 
The functions $h$ and $N$ are, respectively, in  $\mathcal{A}^{q,\gamma_q}$ and  $\mathcal{A}^{p,\gamma_p}$ with  
$p>2,$ $q>\frac{2p}{p-2}$, $\gamma_p>\frac{2d(q+p)}{q}$ and $
\gamma_q>\frac{2d(q+p)}{p};$ 
for every $x$ in $\mathbb{R}^d,$ $N(0,x)=0,$ and $(N(t,x),\interval)$ is a continuous
\filtrationf-local martingale. 
Then
$M(\cdot,L)$ has finite quadratic variation and
$
\pqa M(\cdot,L)\pqc=\pqa M(\cdot,x)\pqc_{x=L}.
$
\end{prop}

\begin{proof}
We have to show that hypotheses of proposition \ref{p5.9} are satisfied.
We denote with $(T_k)_{k\in \mathbb{N}}$ a suitable sequence of $\mathbb{F}$-stopping times such that $N$ and $h$ belong respectively to $\bigcap_{k\in \mathbb{N}}C^{p,\gamma_p}_{T_k}$ and $\bigcap_{k\in \mathbb{N}}C^{q,\gamma_q}_{T_k}.$ Let $(S_k)_{k\in \mathbb{N}}$ be a suitable sequence of $\mathbb{F}$-stopping times such that, for every $k$ in $\mathbb{N},$ $S^k$ is the first instant, between $0$ and $1,$ the process $\vaa N(\cdot,0)\vac+\vaa h(\cdot,0)\vac$ is greater than $k.$ Set, for every $k$ in $ \mathbb{N},$
$
R_k=S_k\wedge T_k,\quad N^k=N^{R_k}$, $M^k=M^{R_k}$ $h^k=h^{R_k}$,
and
$
\Omega_k=\left\{\sup_{ t\in[0,1]}\left|N(t,0)\right|\leq
k \pgc \cap \pga \sup_{t\in[0,1]}\left|h(t,0)\right|\leq k \pgc\cap \pga  R_k=1\pgc.
$
Let $C$ be a compact set of $\mathbb{R}^d,$ $x,y$ in $C,$ and $k$ in $\mathbb{N}$.
Using arguments already employed in the proof of lemma \ref{l5.4}, it is not difficult to show that
if $\alpha=\frac{pq}{p+q}>2,$ there exists a constant $d_k>0,$ depending on $C$ and $k,$ such that,
\beqn \nonu
\mathbb{E}\pqa \sup_{t \in [0,R_k]}\vaa {M}(t,x)-{M}(t,y)\vac^\alpha \pqc&=&\mathbb{E}\pqa \sup_{t \in [0,1]}\vaa {M}^k(t,x)-{M}^k(t,y)\vac^\alpha\pqc \\ \nonu
&\leq&
d_k\vaa x-y\vac^\gamma,
\eeqn
with  $\gamma=\frac{\gamma_p \alpha}{p} \wedge \frac{\gamma_q \alpha}{q}>2d.$ This concludes the proof.
\end{proof}

>From  proposition \ref{p5.9} and lemma \ref{l4.17} we can easily derive  the following corollary.

\begin{coro}\label{c4.25}
Suppose that
$
M(\cdot,x)=\int_0^\cdot h(t,x)dN
$
where $N$ is a continuous $\mathbb{F}$-local martingale and $h$ belongs to $\mathcal{A}^{q,\gamma_q},$ $q>2,$ $\gamma_q>2d.$ 
Then the function $x \mapsto \pqa M(\cdot,x), M(\cdot,x)\pqc$
admits a continuous version, the process
$M(\cdot,L)$ has finite quadratic variation and
$
\pqa M(\cdot,L), M(\cdot,L)\pqc=\pqa M(\cdot,x), M(\cdot,x)\pqc_{x=L}.
$
\end{coro}
\subsubsection{Chain rule formula}
We conclude this section by proving the associative property of forward integrals for the processes studied in this part of the paper.

\begin{prop}\label{c4.26}
Let $h$ and $k$ be respectively in $\mathcal{A}^{p,\gamma_p}$ and  $\mathcal{A}^{q,\gamma_q},$  with  $p>1,$ $q>\frac{p}{p-1}$, $\gamma_p>\frac{d(q+p)}{q},$ $
\gamma_q>\frac{d(q+p)}{p}.$ 
Let
$N$ be  a continuous \filtrationf-local martingale. 
Then $x \mapsto  \int_0^\cdot k(t,x)d^-\int_0^t h(s,x)dN_s=\int_0^\cdot h(t,x)k(t,x)dN_s$ admits a continuous version, $\int_0^\cdot k(t,L)d^-\int_0^t h(s,L)d^-N_s$ exists as limit ucp of its regularizations and
\beqn \nonu
 \int_0^\cdot k(t,L)d^-\int_0^t h(s,L)dN_s&=&\int_0^\cdot k(t,L)h(t,L)d^-N_t. \\ \nonu
\eeqn
\end{prop}
\begin{proof}
By point 3. of lemma \ref{l4.17} we know that $\int_0^\cdot h(t,x)dN_t$ belongs to $\mathcal{A}^{p,\gamma_p}.$ Then, by proposition \ref{p5.5}, $x \mapsto  \int_0^\cdot k(t,x)d\int_0^t h(s,x)dN_s=\int_0^\cdot k(t,x)h(t,x)dN_t$ admits a continuous version and
\beqn \nonu
 \int_0^\cdot k(t,L)d^-\int_0^t h(s,L)dN_s&=&\pta \int_0^\cdot k(t,x)d^-\int_0^t h(s,x)dN_s\ptc_{x=L}\\ \nonu
 &=&\pta \int_0^\cdot k(t,x)h(t,x)dN_t\ptc_{x=L}.
\eeqn
Point 2. of lemma \ref{l4.17} again shows that $hk$ satisfies the hypotheses of corollary \ref{c5.6}.
As a consequence of of this
$$
\int_0^\cdot h(t,L)k(t,L)d^-N_t=\pta \int_0^\cdot k(t,x)h(t,x)dN_t\ptc_{x=L},
$$
and we achieve the end of the proof.
\end{proof}
\begin{defi}\label{d3.33}
\begin{enumerate}
\item
For $p>1$ and $\gamma>d,$ we define $\mathcal{A}^{p,\gamma}(L)$ as the set of all processes $(h(t,L), \interval)$ with $h$ belonging to $\mathcal{A}^{p,\gamma}.$       
\item
$\mathcal{S}(\mathcal{A}^{p,\gamma}(L))$ will be the set of all processes $
h=h^0I_{\{0\}}+\sum_{i=1}^{m}h^iI_{(t_{i-1},t_{i}]}
$
where
$0=t_0\leq t_1,\cdots,t_m=1,$ and $h^i$ belongs to $\mathcal{A}^{p,\gamma}(L)$ for $i=1,...,m.$
\end{enumerate}
\end{defi}

Using similar arguments employed in the proof of proposition \ref{3.8}, it is possible to  demonstrate the following one. 
\begin{prop}\label{p3.35}
Let $N$ be a continuous $\mathbb{F}$-local martingale and $h$ be a process in $\mathcal{S}(\mathcal{A}^{p,\gamma}(L))$ for some $p>1,$ and $\gamma>d.$  Then the following statements are true.
\begin{enumerate}
\item  $h$ is  $N$-forward integrable and the forward integral $\int_0^\cdot h_td^-N_t$ still belongs to $\mathcal{S}(\mathcal{A}^{p,\gamma}(L))$.
\item If $p>2$ and $\gamma>2d,$ $\int_0^\cdot h_td^-N_t$  has finite quadratic variation equal to $\int_0^\cdot h^2_td\pqa N\pqc_t.$
\item If  $k$ belongs to $\mathcal{S}(\mathcal{A}^{q,\gamma_q}(L))$  with  $q>\frac{p}{p-1}$, $\gamma_p>\frac{d(p+q)}{q}$ and $\gamma_q>\frac{d(p+q)}{p},$ then $k$ is forward integrable with respect to $\int_0^\cdot h_td^-N_t$ and 
$$
\int_0^\cdot k_td^-\int_0^t h_sd^-N_s=\int_0^\cdot k_th_td^-N_t.
$$
\end{enumerate}
\end{prop}

\section{$\mathcal{A}$-martingales}
Throughout this section $\cal A$ will be a real linear set of measurable processes indexed by $[0,1)$ with paths which are bounded on each compact interval of $[0,1).$

We will denote with $\mathbb{F}=\pta \mathcal{F}_t\ptc_{t\in[0,1]}$ a filtration indexed by $[0,1]$ and with  $\mathcal{P}(\mathbb{F})$ the $\sigma$-algebra generated by all left continuous and $\mathbb{F}$-adapted processes. In the remainder of the paper we will  adopt the notations $\mathbb{F}$ and $\mathcal{P}(\mathbb{F})$ even when the filtration  $\mathbb{F}$ is indexed by $[0,1).$ At the same way, if $X$ is a process indexed by $[0,1],$ we shall continue to denote with  $X$ its restriction to $[0,1).$ 
\subsection{Definitions and properties}
\begin{defi} \label{d3.1}
A process $X=(X_t,\interval)$ is said
\textbf{$\mathcal{A}$-martingale} if 
every $\theta$ in $\mathcal{A}$
is $X$-improperly forward integrable 
and 
$\mathbb{E}\pqa \int_0^t \theta_sd^-X_s\pqc=0$ for every $\interval.$
\end{defi}

\begin{defi}\label{d4.2}
A process $X=(X_t,\interval)$ is said
\textbf{$\mathcal{A}$-semimartingale} if  it can be written as the sum of an $\cal{A}$-martingale $M$ and a bounded variation process $V,$ with $V_0=0.$
\end{defi}

\begin{rema}\label{r3.2}
\begin{enumerate}
\item
If $X$ is a continuous $\mathcal{A}$-martingale with $X$ belonging to $\cal{A},$ its quadratic variation exists improperly. In fact, if $\int_0^\cdot X_td^-X_t$
exists improperly, it is possible to show that $\pqa X,X\pqc$ exists improperly  and
$
\pqa X,X\pqc=X^2-X_0^2-2\int_0^\cdot X_sd^-X_s.
$
We refer to proposition 4.1 of \cite{RV2} for details.
\item Let $X$ a continuous square integrable martingale with respect to some filtration $\mathbb{F}.$ Suppose that every process in $\cal{A}$ is the restriction to $[0,1)$ of a process $(\theta_t,\interval)$ which is $\mathbb{F}$-adapted, it has  left continuous with right limit paths and $\mathbb{E}\pqa \int_0^1 \theta_t^2 d\pqa X \pqc_t\pqc<+\infty.$ Then $X$ is an $\mathcal{A}$-martingale.  

\item In \cite{FWY} the authors introduced the notion of
  \textbf{weak-martingale}. 
   A semimartingale  $X$ is a
  weak-martingale if $\mathbb{E}\pqa \int_0^t f(s,X_s)dX_s
  \pqc=0,$ $\interval,$ for every $f:\mathbb{R}^+\times \mathbb{R}\rightarrow \mathbb{R},$ bounded
  Borel-measurable. Clearly we can affirm the following.  Suppose that
  $\mathcal{A}$ contains all processes of the form $f(\cdot,X),$ with
  $f$ as above. Let $X$ be a semimartingale $X$ which is an $\mathcal{A}$-martingale. Then $X$ is a weak-martingale. 
\end{enumerate}
\end{rema}

\begin{prop} \label{p3.4a}
Let  $X$ be a continuous $\mathcal{A}$-martingale. The following statements hold true.
\begin{enumerate}
    \item
If  $X$ belongs to $\cal{A},$ $X_0=0$ and $\pqa X,X\pqc=0.$ Then $X\equiv 0
.$
\item Suppose that $\mathcal{A}$ contains all bounded  $\mathcal{P}(\mathbb{F})$-measurable processes.
Then $X$ is an $\mathbb{F}$-martingale.
\end{enumerate}
\end{prop}
\begin{proof}
>From point 1. of remark \ref{r3.2}, $\mathbb{E}\pqa X_t^2\pqc=0,$ for all $\interval.$  

Regarding point 2. it is sufficient to observe that processes of type $I_{A}I_{(s,t]},$ with $0  \leq s \leq t \leq 1,$ and $A$ in $\mathcal{F}_s$ belong to $\cal A.$ Moreover $\int_0^1 I_AI_{(s,t]}(r)d^-X_r=I_A(X_t-X_s).$ This imply $\mathbb E [X_t-X_s\left|\right.\mathcal{F} _s]=0,$ $0\leq s \leq t \leq 1.$ 
\end{proof}

\begin{coro}
The decomposition of an $\cal A$-semimartingale $X$ in definition \ref{d4.2} is unique among the class of processes of type $M+V,$ being $M$ a continuous $\mathcal{A}$-martingale in $\cal A$ and $V$ a bounded variation process.
\end{coro}
\begin{proof}
If $M+V$ and $N+W$ are two decompositions of that type, then $M-N$ is a continuous $\mathcal{A}$-martingale in $\cal A$ starting at zero with zero quadratic variation.  Point 1. of proposition \ref{p3.4a} permits to conclude.
\end{proof}

The following proposition gives sufficient conditions for an $\cal{A}$-martingale to be a martingale with respect to some filtration $\mathbb{F},$ when $\cal{A}$ is made up of $\mathcal{P}(\mathbb{F})$-measurable processes. It constitutes a generalization of point 2. in proposition \ref{p3.4a}.

\begin{defi}\label{d3.6}
We will say that $\mathcal{A}$ satisfies \textbf{assumption $\mathcal{D}$} with respect to a filtration $\mathbb{F}$ if 
\begin{enumerate}
\item Every $\theta$ in $\mathcal{A}$ is $\mathbb{F}$-adapted;
\item For every $0\leq s <1$ there exists a basis  $\mathcal{B}_s$ for $\mathcal{F}_s,$ with the following property. For every $A$ in $\mathcal{B}_s$ there exists a sequence of $\mathcal{F}_s$-measurable random variables $(\Theta_n)_{n\in \mathbb{N}},$ such that for each $n$ the process $\Theta_nI_{[0,1)}$ belongs to $\mathcal{A},$ $\sup_{n\in \mathbb{N}}\vaa \Theta_n\vac\leq 1,$ almost surely and
$$
\lim_{n\rightarrow +\infty}\Theta_n=I_{A}, \quad a.s.
$$
\end{enumerate}
\end{defi}

\begin{prop}\label{p3.3}
Let $X=(X_t,\interval)$ be a continuous {$\mathcal{A}$-martingale} adapted to some filtration $\mathbb{F}$, with $X_t$ belonging to $L^1(\Omega)$ for every $\interval.$ Suppose that $\mathcal{A}$ satisfies assumption $\mathcal{D}$ with respect to $\mathbb{F}.$ Then $X$ is an \filtrationf-martingale.
\end{prop}

\begin{proof}
We have to show that for all $0\leq s\leq t\leq 1$,
$\mathbb{E}\pqa I_{A}\pta X_t-X_s\ptc
\pqc=0$, for all  $A$ in $\mathcal{F}_s$.
We fix $0\leq s<t\leq  1$ and $A$ in $\mathcal{B}_s.$  Let $(\Theta_n)_{n\in \mathbb{N}}$ be a sequence of random variables converging almost surely to $I_{A}$ as in the hypothesis. Since $X$ is an $\cal{A}$-martingale,
$\mathbb{E}^{}\pqa \Theta_n \pta X_t-X_s\ptc
\pqc=0$, for all $n$ in $\mathbb{N}$.
We note that $X_t-X_s$ belongs to $L^1(\Omega),$ then, by Lebesgue dominated convergence theorem,
\beqn\nonu
\vaa \mathbb{E}\pqa I_{A}\pta
X_t-X_s\ptc\pqc\vac\leq\lim_{n\rightarrow +\infty}\mathbb{E}\pqa
\vaa I_{A}-\Theta_n \vac \vaa X_t-X_s \vac\pqc=0.
\eeqn
Previous result extends to the whole $\sigma$-algebra $\mathcal{F}_s$ and this permits to achieve the end of the proof.  
\end{proof}

Some interesting properties can be derived taking inspiration from
\cite{FWY}.

For a process $X,$ we will denote
\beqn \label{AX}
\mathcal{A}_X&=&\pga (\psi(t,X_t)), \intervala
\left |\right. \psi:[0,1]\times \mathbb{R} \rightarrow \mathbb{R}, \mbox{ Borel-measurable }\right. \\ \nonu &&\left.
\mbox{ with polynomial growth and lower bounded}\pgc.
\eeqn
\begin{rema}
At this stage we could avoid to impose a lower bound on functions in $\mathcal{A}_X.$ Nevertheless, we prefer to consider this qualitative restriction in view of further applications to finance. Indeed, $\mathcal{A}_X$ will play the role of a possible class of \textit{admissible portfolios} and we are interested in excluding among them the so called \textit{doubling strategies}. Generally speaking, a doubling strategy is an \textit{arbitrage} which can be realized if unbounded accumulation of losses are allowed. For more details about this arguments the reader is referred to Harrison and Pliska (1979).
  
\end{rema}
\begin{prop}
Let $X$ be a continuous $\mathcal{A}$-martingale with $\mathcal{A}=\mathcal{A}_X.$ 

Then, for every $\psi$ in $C^2(\mathbb R)$ with bounded first and second derivatives, the process 
$$
\psi(X)-\frac{1}{2}\int^{\cdot}_{0}\psi''(X_s)d\pqa X,X\pqc_s
$$
is an $\mathcal{A}$-martingale.
\end{prop}
\begin{proof}
The process $X$ belongs to $\cal A.$ In particular, $X$ admits improper quadratic variation. We set $Y=\psi(X)-\frac{1}{2}\int^{\cdot}_{0}\psi''(X_s)d\pqa X,X\pqc_s.$
Let $\theta$ in $\mathcal {A}_X.$ By lemma \ref{l2.1}, for every $0\leq t < 1$
$$
\int^{t}_{0}\theta_s d^-Y_s =\int^{t}_{0}\theta_s\psi'(X_s)d^-X_s.
$$
Since $\theta \psi'(X)$ still belongs to $\cal A,$ $\theta$ is $Y$-improperly integrable and  
\beqn \label{e8}
\int^{\cdot}_{0}\theta_t d^-Y_t =\int^{\cdot}_{0}\theta_t\psi'(X_t)d^-X_t.
\eeqn
We  conclude taking the expectation in equality (\ref{e8}). 
\end{proof}

\begin{prop}\label{p4.9}
Suppose that $\cal A$ is an algebra. Let $X$ and $Y$ be two continuous $\mathcal{A}$-martingales with $X$ and $Y$ in $\cal{A}.$ 

Then the process $XY-\pqa X,Y\pqc$ is an $\mathcal{A}$-martingale . 
\end{prop}

\begin{proof}
Since $\cal A$ is a real linear space, $(X+Y)$ belongs to $\cal A.$ In particular by point 1. of remark \ref{r3.2}, $\pqa X+Y,X+Y\pqc,$ $\pqa X,X\pqc$ and $\pqa Y,Y\pqc$ exist improperly. 
This implies that $\pqa X,Y\pqc$ exists improperly too and that it is a bounded variation process. 
Therefore the vector $(X,Y)$ admits all its mutual brackets on each compact set of $[0,1).$ 
Let $\theta$ be in $\mathcal{A}.$ 
Since
$\mathcal{A}$ is an algebra, $\theta X$ and $\theta Y$ belong to $\mathcal{A}$ and so both $\int_0^\cdot  \theta_s X_s d^-Y_s$ and $\int_0^\cdot \theta_s Y_s d^-X_s$ locally exist.  
By lemma \ref{l2.1} $\int_0^\cdot \theta_t d^-\pta X_tY_t-\pqa X,Y\pqc\ptc $ exists improperly too  and
\beqn \nonu
\int_0^\cdot \theta_td^-\pta X_tY_t-\pqa X,Y\pqc_t\ptc=\int_0^\cdot Y_t\theta_t d^-X_t+\int_0^\cdot  X_t\theta_t  d^-Y_t. 
\eeqn
Taking the expectation in the last expression we then get the result.
\end{proof}

We recall a notion and a related result of \cite{CR}. 

A process $R$ is {\textit{{strongly predictable}}} with respect to a filtration $\mathbb{F},$ if 
$$
\exists\ \delta>0, \mbox{ such that } R_{\ep+\cdot} \mbox{ is } \mathbb{F}\mbox{-adapted},\ \mbox{ for every } \ep\leq\delta.
$$

\begin{prop} \label{p3.11}
Let $R$ be an $\mathbb{F}$-strongly predictable continuous process. Then for every continuous $\mathbb{F}$-local martingale $Y,$ $\pqa R,Y\pqc=0$. 
\end{prop}

Proposition \ref{p3.11} combined with  proposition \ref{p4.9} implies proposition \ref{p3.10} and corollary \ref{c3.11}.

\begin{prop}\label{p3.10}
Let $\cal{A},$ $X$ and $Y$ be as in proposition \ref{p4.9}. Assume, moreover, that $X$ is an $\mathbb{F}$-local martingale, and that  $Y$ is strongly predictable with respect to $\mathbb{F}.$
Then $XY$ is an $\mathcal{A}$-martingale.
\end{prop}

\begin{coro}\label{c3.11}
Let $\cal{A},$ $X$ and $Y$ be as in proposition \ref{p4.9}.
Assume that $X$ is a local martingale with respect to some
filtration $\mathbb{G}$ and that $Y$ is either $\mathbb{F}$-independent, or $\mathcal{G}_0$-measurable. 
Then $XY$ is an $\mathcal{A}$-martingale.
\end{coro}
\begin{proof}
If $Y$ is $\mathbb{G}$-independent, it is sufficient to apply previous proposition with
$\mathbb{F}=\pta \bigcap_{\ep>0} \mathcal{G}_{t+\ep}\vee \sigma(Y)\ptc_{t\in[0,1]}.
$
\end{proof}

\subsection{$\mathcal{A}$-martingales and Weak Brownian motion}
We proceed defining and discussing processes which are \textit{weak-Brownian motions} in order to exhibit explicit examples of $\mathcal{A}$-martingales.

\begin{defi}(\cite{FWY})
A stochastic process $(X_t, \interval)$ is a \textbf{weak Brownian
  motion of order $k$} if for every k-tuple $(t_1,t_2,...,t_k)$
$$
(X_{t_1},X_{t_2},\dots,X_{t_k})\stackrel{law}{=}(W_{t_1},W_{t_2},\dots,W_{t_k})
$$
where $(W_t,\interval)$ is a Brownian motion.
\end{defi}
We set, for a process $(X_t, \interval),$
\beqn \nonu
\mathcal{A}^1_{X}&=&\pga (\psi(t,x), \interval, \mbox{ with polynomial growth  s.t } \psi=\partial_x \Psi \right. \\ \nonu && \left.\Psi \in C^{1,2}([0,1]\times \mathbb{R}) \mbox{ with  }  \vaa \partial_t\Psi\vac  +\vaa \partial_{xx}\Psi\vac \mbox{ bounded }\pgc.
\eeqn
\begin{assu}\label{a4.15}
We suppose that $\sigma:[0,1]\times
\mathbb{R}\rightarrow \mathbb{R}$ is a 
Borel-measurable and bounded function such that the following equation has a unique solution $\pta \nu_t\ptc_{t\in [0,1]}$ in the sense of distribution
\beqn \label{e11}
\left\{
\begin{array}{ll}
\partial_t \nu_t(dx)=\frac{1}{2}\pta \sigma(t,x)^2 \nu_t(dx)\ptc^{''}\\
\nu_0(dx)=\delta_{0}.
\end{array} \right.
\eeqn
\end{assu}

\begin{rema}
Assumption \ref{a4.15} is verified for $\sigma(t,x)\equiv \sigma,$ being  $\sigma$ a positive real constant and,  in that case, $\nu_t=N(0,\sigma^2t),$ for every $\interval.$   
\end{rema}

\begin{prop}\label{p}
Let $(X_t,\interval)$ be a continuous finite quadratic variation process with $X_0=0,$  and
$d\pqa X\pqc_t=(\sigma(t,X_t))^2 dt,$ where $\sigma$ fulfills assumption \ref{a4.15}.  
Then the following statements are true.
\begin{enumerate}
\item \label{point1}
Suppose that $\mathcal A=\mathcal{A}^1_X.$
Then  $X$ is an $\mathcal{A}$-martingale if and only if, for every 
  $\interval,$ $X_t\stackrel{law}{=}Z_t,$ for every $(Z,B)$ solution of
  equation
$
dZ=\sigma(\cdot,Z)dB, Z_0=0,
$
in the sense of definition 1.2 in chapter $IX$ of \cite{RY}. In particular, if
$\sigma\equiv 1,$ $X$ is a weak Brownian motion of order $1$, if and only if it is an $\mathcal{A}_X^{1}$-martingale.
\item \label{point2}
Suppose that $d\pqa X\pqc_t=f_tdt,$ with $f$ $\mathcal{B}([0,1])$-measurable and bounded. If  $X$ is a weak Brownian motion
of order $k=1,$ then $X$ is an $\mathcal{A}$-semimartingale. Moreover the process 
$$
X+\int_0^\cdot \frac{(1-f_s)X_s}{2s} ds.
$$
is an $\mathcal{A}$-martingale.
\end{enumerate}
\end{prop}

\begin{proof}
\begin{enumerate}
\item
Using It\^o inverse
formula recalled in proposition \ref{p2.1} we can write, for
every $\interval$ and  $\psi=\partial_x\Psi$ in $\mathcal{A}^1_X$
\begin{eqnarray}\nonu
\int_0^t \psi(s,X_s)d^-X_s&=&\Psi(t,X_t)-\Psi(0,X_0)  \\ &-&  \int_0^t \pta \partial_s\Psi
+\frac{1}{2}\partial_{xx}^{(2)}\Psi\sigma^2\ptc(s,X_s) ds. \label{f2}
\end{eqnarray}
For every $\interval,$ we denote with $\mu_t(dx)$ the law of $X_t.$ If $X$ is an $\mathcal{A}^1_X$-martingale, from (\ref{f2}) we derive      
\beqn \nonu
0&=&\int_\mathbb{R} \Psi(t,x)d\mu_t(x)-\int_\mathbb{R}  \Psi(0,x)\mu_0(dx) -\int_0^t\int_\mathbb{R}
\partial_s \Psi(s,x)\mu_s(dx) ds\\ \label{e3} &-&\frac{1}{2}\int_0^t \int_\mathbb{R}
\partial_{xx}^{(2)}\Psi(s,x)\sigma(s,x)^2 \mu_s(dx) ds.
\eeqn

In particular,  the law of $X$ solves equation (\ref{e11}).

On the other hand, let
$(Z,B)$ be a solution of equation
$
Z=\int_0^\cdot \sigma(s,Z_s)dB_s.
$ 
The law of $Z$ fulfills equation (\ref{e3}) too. Indeed, $Z$ is a finite quadratic variation process with $d\pqa Z\pqc_t=(\sigma(t,Z_t))^2dt$ which is an $\mathcal{A}^1_X$-martingale by point 2. of remark \ref{r3.2}.
 By assumption (\ref{a4.15}) $X_t$ must have the same law of $Z_t.$ This establishes the converse implication of point 1.

Suppose, on the contrary, that $X_t$ has the same law of $Z_t,$ for every $\interval.$ Using the fact that $Z$ is an $\mathcal{A}^1_X$-martingale which solves equation (\ref{f2}) we get
$$
\mathbb{E}\pqa \Psi(t,Z_t)-\Psi(0,Z_0)-\int_0^t \pta \partial_s\Psi 
+\frac{1}{2}\partial_{xx}^{(2)}\Psi \sigma^2\ptc (s,X_s) ds \pqc=0.
$$
for every $\Psi$ in $C^{1,2}([0,1]\times \mathbb{R})$ with $\partial_x \Psi=\psi$ in $\mathcal{A}^1_X.$
Since $X_t$ has the same law of $Z_t,$ for every $\interval,$ equality (\ref{f2}) implies that  $$
\mathbb{E}\pqa \int_0^\cdot \psi(t,X_t)d^-X_t\pqc=\mathbb{E}\pqa \int_0^\cdot \psi(t,Z_t)d^-Z_t\pqc=0, 
$$
The proof of the first point is now  achieved.
\item
Suppose that $\sigma(t,x)^2=f_t,$ for every $(t,x)$ in $[0,1]\times
\mathbb{R}.$ Let $\Psi$ be in $C^{1,2}\pta [0,1]\times \mathbb{R}\ptc$ such that $\psi=\partial_x \Psi$ belongs to $\mathcal{A}^1_X.$ Proposition \ref{p2.1} yields 
$$
\int_0^t \psi(s,X_s)d^-X_s=Y^\Psi_t+\frac{1}{2}\int_0^t
\partial_{xx}^{(2)}\Psi(s,X_s)(1-f_s)ds, \quad \interval,
$$
with
$$
Y^\Psi_t=\Psi(t,X_t)-\Psi(0,X_0)-\int_0^t
\partial_{s}\Psi_s(s,X_s)ds-\frac{1}{2}\int_0^t \partial_{xx}^{(2)}\Psi(s,X_s)ds.
$$
Moreover $X$ is a weak Brownian motion of order
$1.$ This implies 
$
\mathbb{E}\pqa Y^\Psi_t \pqc=0$, for every $\interval$.
We derive that 
$$
\mathbb{E}\pqa \int_0^t \psi(s,X_s)d^-X_s+\frac{1}{2}\int_0^t
\partial_{xx}^{(2)}\Psi(s,X_s)(f_s-1)ds\pqc=\mathbb{E}\pqa Y^\Psi_t \pqc=0.
$$
Since the law of $X_t$ is $N(0,t),$ by Fubini theorem and integration by parts on the real line we obtain 
$$
\mathbb{E}\pqa \int_0^t
\partial_{xx}^{(2)}\Psi(s,X_s)(f_s-1)ds\pqc=\mathbb{E}\pqa \int_0^t
\psi(s,X_s)\frac{(1-f_s)X_s}{s}ds\pqc.
$$
This concludes the proof of the second point.
\end{enumerate}
\end{proof}

>From \cite{FWY} we can extract an example of an $\cal{A}$-semimartingale  which is not a semimartingale.

\begin{exam}
Suppose that $\pta B_t, \interval\ptc$ is a Brownian motion on the probability space $\pta \Omega, \mathbb{G},P\ptc,$ being $\mathbb{G}$ some filtration on $(\Omega, \mathcal{F},P).$ Set
$$
X_t=\pga
\begin{array}{ll}
B_t,& 0\leq t\leq \frac{1}{2} \\
B_{\frac{1}{2}}+(\sqrt{2}-1)B_{t-\frac{1}{2}}, & \frac{1}{2}< t\leq 1.
\end{array}
\right.
$$
Then $X$ is a continuous weak Brownian motion of order $1,$ which is not a $\mathbb{G}$-semimartingale. Moreover it is possible to show that $d\pqa X \pqc_t=f_tdt,$ with $f=I_{[0,\frac{1}{2}]}+(\sqrt{2}-1)^2I_{[\frac{1}{2},1]}.$
In particular, thanks to point 2. of previous proposition, $X+\int_0^\cdot \frac{(1-f_s)X_s}{2s}ds$  is an ${\cal{A}}^1_X$-martingale. 
\end{exam}

A natural question is the following. Supposing that $X$ is an $\cal{A}$-martingale with respect to a probability measure $Q$ equivalent to $P,$ what can we say about the nature of $X$ under $P$ ?. The following proposition provides a partial answer to this problem when $\mathcal{A}=\mathcal{A}_X^1.$
  
\begin{prop} 
Let $X$ be as in proposition \ref{p}, and $\sigma$  satisfy assumption \ref{a4.15}. Assume, furthermore,  that  $X$ is an  $\mathcal{A}^1_X$-martingale under a
probability measure $Q$ with $P << Q,$   
Then the law of $X_t$ is absolutely continuous with respect to
Lebesgue measure, for every $\interval.$
\end{prop}
\begin{proof}
Since $P << Q,$ for every $\interval,$ $X_t(P)<<X_t(Q).$ Then  it is sufficient to observe that by proposition \ref{p}, for every
  $\interval,$ the law of $X_t$ under $Q$ is absolutely continuous with respect to Lebesgue. 
\end{proof}
  
\begin{coro} \label{c4.16}
Let $X$ be as in proposition \ref{p}, and $\sigma$  satisfy assumption \ref{a4.15}. Assume, furthermore,  that  $X$ is an  $\mathcal{A}_X$-martingale under a
probability measure $Q$ with $P << Q,$   
Then the law of $X_t$ is absolutely continuous with respect to
Lebesgue measure, for every $\interval.$
\end{coro}  
\begin{proof}
Clearly $\mathcal{A}^1_X$ is contained in $\mathcal{A}_X.$ The result is then a consequence of previous proposition. 
\end{proof}
\begin{prop}
Let $\pta X_t,\interval\ptc$ be a continuous weak Brownian motion of
order $ 8.$ Then, for every $\psi:[0,1]\times \mathbb{R} \rightarrow \mathbb{R}$, Borel measurable with polynomial growth,  
the forward integral $\int_0^\cdot \psi(t,X_t)d^-X_t,$ exists and
$$
\mathbb{E}\pqa \int_0^\cdot\psi(t,X_t)d^-X_t, \pqc=0.
$$
In particular, $X$ is an $\mathcal{A}_X$-martingale.
\end{prop}

\begin{proof}
Let $\psi:[0,1]\times \mathbb{R} \rightarrow \mathbb{R}$ be Borel measurable and  $t$ in $\interval$ be fixed. Set
$$
I_\ep^X(t)=I(\ep,\psi(\cdot,X),X)
\quad I_\ep^B(t)=I(\ep,\psi(\cdot,B),B),
$$
being $B$ a  Brownian motion on a filtered probability space $(\Omega^B, \mathbb{F}^B,P^B).$ 

Since $X$ is a weak Brownian motion of order $8,$ it follows that  
$$
\mathbb{E}\pqa \vaa I_\ep^X(t)-I_\delta^X(t)\vac^4\pqc=\mathbb{E}^{P^B}\pqa \vaa I_\ep^B(t)-I_\delta^B(t)\vac^4\pqc, \quad  \forall \ \ep, \delta>0.
$$
We show now that $I^B_\ep(t)$ converges in $L^4(\Omega).$ This implies that $I^X_\ep(t)$ is of Cauchy in $L^4(\Omega).$

In \cite{RV05}, chapter 3.5, it is proved that $I_\ep^B(t)$ converges in probability when $\ep$ goes to zero, and the limit equals the It\^o integral $\int_0^t \psi(s,B_s)dB_s.$
Applying Fubini theorem for It\^o integrals, theorem 45 of \cite{P}, chapter IV and  Burkholder-Davies-Gundy inequality, 
we can perform the following estimate, for every $p>4:$
\beqn \nonu
\mathbb{E}^{P^B}\pqa \vaa I_\ep^B(t)\vac^p\pqc&=&\mathbb{E}^{P^B}\pqa \vaa \int_0^t \pta \frac{1}{\ep}\int_{r-\ep}^r \psi(s,B_s)ds \ptc dB_r \vac^p\pqc \\ \nonu
&\leq& c\mathbb{E}^{P^B} \pqa  \int_0^1 \frac{1}{\ep} \int_{r-\ep}^r \vaa \psi(s,B_s)\vac^p ds dr\pqc \\ \nonu 
&\leq& c \sup_{t\in[0,1]}\mathbb{E}^{P^B}\pqa \vaa \psi(t,B_t)\vac^p\pqc< +\infty,
\eeqn
for some positive constant $c.$
This implies the uniformly integrability of the family of random variables $\pta  (I_\ep^B(t))^4 \ptc_{\ep>0}$ and therefore the convergence in $L^4(\Omega^B,P^B)$ of $\pta I_\ep^B(t)\ptc_{\ep>0}.$

Consequently,  $\pta  I_\ep^X(t) \ptc_{\ep>0}$ converges in $L^4(\Omega)$ toward a random variable $I(t).$ It is clear that $\mathbb{E}\pqa I(t)\pqc=0,$ being $I(t)$ the limit in $L^2(\Omega)$ of random variables having zero expectation. 

To conclude we show that Kolmogorov lemma applies to find a continuous version of $\pta I(t), \interval\ptc.$ Let $0\leq s\leq t\leq 1.$
Applying the same arguments used above
\beqn \nonu
\mathbb{E}\pqa \vaa I(t)-I(s)\vac^4\pqc
&=&\lim_{\ep \rightarrow
  0}\mathbb{E}^{P^B} \pqa \vaa \int_s^t \pta \frac{1}{\ep}\int_{u}^{u+\ep}\psi(u,B_u)dB_r \ptc du \vac^4\pqc \\ \nonu
  &\leq& c \mathbb{E}^{P^B} \pqa \vaa \int_s^t \pta\frac{1}{\ep}\int_{r-\ep}^r \psi(u,B_u)du \ptc^2 dr \vac^2\pqc \\ \nonu
 &\leq& \vaa t-s \vac \mathbb{E}^{P^B} \pqa \int_s^t \frac{1}{\ep}\int_{r-\ep}^r \vaa \psi(u,B_u)\vac^4 du  dr  \pqc \\ \nonu
&\leq& \sup_{u\in[0,1]} \mathbb{E}^{P^B}\pqa  \vaa \psi(u,B_u)\vac ^4 \pqc \vaa t-s \vac^2, \quad c>0.
\eeqn
\end{proof}

\subsection{Optimization problems and  $\cal{A}$-martingale
  property}
\subsubsection{G\^ateaux-derivative: recalls}
In this part of the paper we recall the notion of G\^ateaux differentiability and we list
some related properties.
\begin{defi}
A function $f:\cal{A}\rightarrow \mathbb{R}$ is said
\textbf{{G\^ateaux-differentiable}} at $\pi \in \cal{A},$ if there
exists $Df_{\pi}:\cal{A}\rightarrow \mathbb{R}$ such that
$$
\lim_{\ep \rightarrow
  0}\frac{f(\pi+\ep\theta)-f(\pi)}{\ep}=Df_{\pi}(\theta),\quad \forall
\theta \in \cal{A}.
$$
If $f$ is \textit{G\^ateaux}-differentiable at every $\pi \in
\cal{A},$ then $f$ is said \textit{G\^ateaux}-differentiable on
$\cal{A}$.
\end{defi}

\begin{defi} \label{d3.18}
Let  $f:\cal{A}\rightarrow \mathbb{R}.$ A process $\pi$ is said
\textbf{optimal} for $f$ in $\mathcal{A}$ if
$$
f(\pi)\geq f(\theta), \quad \forall \theta \in \mathcal{A}.
$$
\end{defi}

We state this useful lemma omitting its straightforward proof.

\begin{lemm}\label{l7.2}
Let  $f:\cal{A}\rightarrow \mathbb{R}.$ For every $\pi$ and $\theta$
in $\cal{A}$ define $f_{\pi,\theta}:
\mathbb{R}\longrightarrow\mathbb{R}$ in the following way:
$$
f_{\pi,\theta}(\lambda)=f(\pi+\lambda(\theta-\pi)).
$$
Then it
holds:
\begin{enumerate}
\item f is \textit{G\^ateaux}-differentiable if and only if for every $\pi$ and $\theta$ in $\cal{A},$  $f_{\pi,\theta}$ is differentiable
  on $\mathbb{R}.$ Moreover
  $f_{\pi,\theta}'(\lambda)=Df_{\pi+\lambda(\theta-\pi)}(\theta-\pi).$
\item  f is concave if and only if  $f_{\pi,\theta}$ is concave
  for every  $\pi$ and $\theta$ in $\cal{A}$.
\end{enumerate}
\end{lemm}

\begin{prop}\label{p7.3}
Let  $f:\cal{A}\rightarrow \mathbb{R}$ be
\textit{G\^ateaux}-differentiable.  Then, if $\pi$ is optimal for $f$
in $\mathcal{A},$ $Df_\pi=0.$ If $f$ is concave
$$
\pi \mbox{ is optimal for } f \mbox{ in }\mathcal{A}  \quad
\Longleftrightarrow  \quad Df_\pi =0.
$$
\end{prop}

\begin{proof}
It is immediate to prove that $\pi$ is optimal for $f$ if and only if $\lambda=0$ is a maximum for $f_{\pi,\theta},$ for every $\theta$ in $\cal{A}.$ By lemma \ref{l7.2} $f^{'}_{\pi,\theta}(0)=Df_{\pi}(\theta),$ for every $\theta$ in $\cal{A}$. The conclusion follows easily.
\end{proof}
\subsubsection{An optimization problem}\label{s3.3.2}
 In this part of the paper $F$ will be supposed to be a measurable function on $(\Omega\times \mathbb{R},\mathcal{F}\otimes \mathcal{B}(\mathbb{R})),$ almost surely in $C^1(\mathbb{R}),$ strictly increasing, with $F'$ being the derivative of $F$ with respect to $x,$ bounded on $\mathbb{R}$, uniformly in $\Omega.$
In the sequel $\xi$ will be a continuous  finite quadratic variation process with $\xi_0=0.$

The starting point of our construction is the following hypothesis.

\begin{assu}\label{a0}
\begin{enumerate}
\item If $\theta$ belongs to $\cal{A},$ then $ \theta I_{[0,t]}$ belongs to $\mathcal{A}$ for every $\intervala.$
\item  Every $\theta$ in $\cal{A}$ $\xi$-improperly forward integrable, and
$$
\mathbb{E}\pqa \vaa \int_0^1 \theta_t d^-\xi_t\vac +\vaa \int_0^1 \theta_t^2 d[\xi]_t \vac \pqc < +\infty.
$$
\end{enumerate}
\end{assu}
\begin{defi} 
Let $\theta$ be  in $\cal{A}.$ We denote 
$$
L^\theta=\int_0^1 \theta_td^-\xi_t-\frac{1}{2}\int_0^1 \theta^2 d[\xi]_t,
\quad 
dQ^\theta=\frac{F'(L^\theta)}{\mathbb{E}\pqa F'(L^\theta)\pqc}
$$ and 
we set $f(\theta)=\mathbb{E}\pqa F(L^\theta) \pqc.$ 
\end{defi}

We observe that point 2. of assumption \ref{a0} and the boundedness of $F'$ implies that  $\mathbb{E}\pqa \vaa F(L^\theta) \vac\pqc <+\infty.$ Therefore $f$ is well defined.
 \begin{rema} Point 2. of assumption \ref{a0} implies that $\mathbb{E}\pqa \vaa \xi_t \vac+\pqa \xi \pqc_t\pqc \le +\infty,$ for every $\interval.$ This is due to the fact that $\cal A$ must contain real constants. 
 \end{rema}
We are interested in describing a link between the existence of an optimal process for $f$ in $\cal{A}$ and the $\mathcal{A}$-semimartingale property for $\xi$ under some probability measure equivalent to $P,$ depending on the optimal process.

\begin{lemm}\label{l5.1}
The function $f$ is \textit{G\^ateaux}-differentiable on
$\cal{A}.$ Moreover for every $\pi$ and $\theta$ in
$\mathcal{A}$
$$
Df_\pi(\theta)=\mathbb{E}\pqa F'(L^\pi)\int_0^1
\theta_t d^-\pta \xi_t- \int_0^t \pi_s d\pqa \xi\pqc_s \ptc \pqc.
$$
If $F$ is concave, then $f$ inherits the property.
\end{lemm}

\begin{proof}
Regarding the concavity of $f,$ we recall that if  $F$ is increasing and concave, it is sufficient to verify that, for every $\theta$ and $\pi$ in $\mathcal{A}$, it holds
$$
L^{\pi+\lambda(\theta-\pi)}-L^{\pi}-\lambda\pta L^{\theta}-L^{\pi}\ptc\geq 0, \quad 0\leq \lambda \leq 1.
$$
A short calculus shows that, for every $0\leq \lambda \leq 1,$
$$
L^{\pi+\lambda(\theta-\pi)}-L^{\pi}-\lambda\pta L^{\theta}-L^{\pi}\ptc=\frac{1}{2} \lambda(1-\lambda)\int_0^1 (\theta_t-\pi_t)^2d\pqa\xi \pqc_t\geq 0.
$$
Using the
differentiability of $F$ we can write
$$
a_\ep=\frac{1}{\ep}(f(\pi+\ep\theta)-f(\pi))=\mathbb{E}\pqa
H^\ep_{\pi,\theta}\int_0^1 F'\pta
L^{\pi}+\mu \ep H^\ep_{\pi,\theta} \ptc
d\mu  \pqc,
$$
with
\beqn \nonu
H^{\ep}_{\pi,\theta}
=\int_0^1 \theta_td^-\xi_t-\frac{1}{2}\int_0^1
(\theta_t^2\ep+2\theta_t\pi_t)d\pqa \xi\pqc_t.
\eeqn
The conclusion follows by Lebesgue dominated convergence theorem, which applies thanks to the boundedness of $F'$ and point 2. in assumption \ref{a0}.
\end{proof}

Putting together lemma \ref{l5.1} and proposition \ref{p7.3} we can
formulate the following.

\begin{prop}\label{p7.8}
If a process $\pi$ in $\cal{A}$ is optimal for  
$
\theta\mapsto \mathbb{E}\pqa F\pta 
L^\theta \ptc \pqc,
$ then  the process
$
\xi-\int_0^\cdot \pi_t d \pqa \xi\pqc_t
$ is an  $\cal{A}$-martingale under $Q^\pi.$ 
If $F$ is concave the converse holds.
\end{prop}

\begin{proof}
Thanks to lemma \ref{l5.1} and point 1. in  assumption \ref{a0}, for every $\theta$ in $\mathcal{A}$ and $\interval$
\beqn \nonu
0&=&Df_\pi(\theta I_{[0,t]})=\mathbb{E}\pqa F'(L^\pi)\int_0^t
\theta_s d^-\pta \xi_s-\int_0^s \pi_r d\pqa \xi\pqc_r \ptc \pqc \\ \nonu
&=&\mathbb{E}^{Q^\pi}\pqa \int_0^t
\theta_s d^-\pta \xi_s-\int_0^s \pi_r d\pqa \xi\pqc_r \ptc \pqc.
\eeqn
\end{proof}

The following proposition describes some sufficient conditions to recover the semimartingale property for $\xi$ with respect to a filtration $\mathbb{G}$ on $(\Omega,\mathcal{F}), $ when the set $\mathcal{A}$ is made up of $\mathbb{G}$-adapted processes. It can be proved using proposition \ref{p3.3}.

\begin{prop}\label{p5.1}
Assume that $\xi$ is adapted with respect to some filtration $\mathbb{G}$ and  that $\mathcal{A}$ satisfies the hypothesis $\mathcal{D}$ with respect to $\mathbb{G}.$
If a process $\pi$ in $\cal{A}$ is optimal for $\theta\mapsto \mathbb{E}\pqa F(L^\theta) \pqc,$ then the process 
$
\xi-\int_0^\cdot \beta_t d \pqa \xi \pqc_t
$
is a $\mathbb{G}$-martingale under $P,$
where   $\beta=\pi+\frac{1}{p^\pi}\frac{d\pqa p^\pi, \xi\pqc}{d\pqa \xi,\xi\pqc},$ and $p^\pi=\mathbb{E}\pqa \frac{dP}{dQ^\pi}\left|\right. \cal{G}_\cdot \pqc.$
If $F$ is concave, then the converse holds. 
\end{prop}
\begin{proof} 
Thanks to point 2. of assumption \ref{a0}, for every $\intervala,$
  the random variable $\xi_t-\int_0^t \pi_t d\pqa\xi \pqc_t$ is in $L^1\pta \Omega \ptc$  and so in $L^1\pta \Omega,Q^\pi\ptc$ being $\frac{dQ^\pi}{dP}$ bounded.  Then proposition \ref{p3.3} applies to state that $\xi-\int_0^\cdot \pi_t d \pqa \xi\pqc_t$ is a $\mathbb{G}$-martingale under $Q^\pi.$  Using Girsanov theorem, chapter 6 of \cite{P}, we get the necessity condition.  As far as the converse is concerned, we observe that, thanks to the hypotheses on $\mathcal{A},$ if $\xi-\int_0^\cdot \pi_t d \pqa \xi\pqc_t
$ is a $\mathbb{G}$-martingale, then for every $\theta$ in $\mathcal{A},$ the process $\int_0^\cdot \theta_t d^-\pta \xi_t- \int_0^t \pi_s d \pqa \xi\pqc_s\ptc$ is a $\mathbb{G}$-martingale starting at zero with zero expectation. This concludes the proof. 
\end{proof}

\begin{prop}\label{c3.18}
Suppose that there exists a measurable process $(\gamma_t,\interval)$ such that the process 
$
\xi-\int_0^\cdot \gamma_t
  d\pqa \xi \pqc_t
$
is an $\cal{A}$-martingale. Assume, furthermore, the existence of a sequence of processes $\pta \theta^n\ptc_{n\in \mathbb{N}}\subset \cal A$ with 
$$
\lim_{n\rightarrow +\infty}\mathbb{E}\pqa \int_0^1 \vaa \theta^n_t-\gamma_t\vac^2d\pqa \xi\pqc_t\pqc=0.
$$ If $\gamma$ belongs to $\cal{A}$ then $\gamma$ is optimal for 
$
\theta\mapsto \mathbb{E}\pqa L^\theta \pqc$. Moreover if there exists an optimal  process $\pi,$ then $d \vaa \pqa \xi \pqc \vac \pga t\in[0,1), \gamma_t\neq \pi_t \pgc=0,$ almost surely.
\end{prop}

\begin{proof}
Using proposition \ref{p7.8} we deduce that a
process $\pi$ is optimal for $f$ if and only if the process
$
\int_0^\cdot (\gamma_t-\pi_t) d\pqa \xi\pqc_t
$
is an  $\cal{A}$-martingale under $P.$
Then $\pi$ is optimal if and only if for every $\theta$ is in $\cal{A}$ it holds: 
$\mathbb{E}\pqa \int_0^\cdot \theta_t (\gamma_t-\pi_t)d\pqa \xi\pqc_t\pqc=0.
$
This permits to achieve immediately the end of the proof.
\end{proof}

\subsubsection{An example of $\cal{A}$-martingale and  a related optimization problem} \label{3.26}
We illustrate a setting where proposition \ref{c3.18} applies. It will be deduced  by \cite{LNN}. In that paper the authors study a particular case
of the optimization problem considered in proposition \ref{c3.18}. As process $\xi$ they take a Brownian motion $W,$ and
they find sufficient conditions in order to have existence of a
process $\gamma$ such that $W-\int_0^\cdot \gamma_tdt$ is an $\cal{A}$-martingale, being $\cal{A}$ some specific set we shall clarify later. To get their goal, they consider an anticipating
setting and combine Malliavin calculus with substitution formulae, the anticipation being generated by a random variable possibly depending on the whole trajectory of $W.$

We work into the specific framework of subsection \ref{s4.1}. 

\begin{assu}\label{a3.21} We suppose the existence of a random variable $L$ in
$\mathbb{D}^{1,2},$ satisfying the following assumption:
\begin{enumerate}
\item $
\int_{\mathbb{R}}\mathbb{E}\pqa \vaa L\vac^2 I_{\pga
  0\leq x\leq L\pgc \cup \pga 0\geq x \geq L \pgc} \pqc dx<+\infty;
$
\item for a.a. $t$ in $[0,1]$ the process
$$
I(\cdot,t,L):=I_{[t,1]}(\cdot)I_{\pga \int_t^1 (D_sL)^2ds>0\pgc}\pta
\int_t^1 (D_sL)^2 ds \ptc^{-1}(D_tL)(D_{\cdot}L)
$$
belongs to $Dom\delta$ and there exists a $\mathcal{P}(\mathbb{F}) \times\mathcal{B}(\mathbb{R})$-measurable random field $\pta h(t,x), \interval, x\in \mathbb{R}\ptc$ such that $h(\cdot,L)$ belongs to $L^2\pta \Omega \times [0,1]\ptc$ and
$$
\mathbb{E}\pqa \int_0^1 I(u,t,L)dW_u \left|\right.\mathcal{F}_t \vee  \sigma(L)\pqc=h(t,L), \quad \interval.
$$
\end{enumerate}
\end{assu}
Let $\Theta(L)$ be the set of processes $\pta \theta_t, \intervala \ptc$ such that there exists a random field $(u(t,x), \interval, x\in \mathbb{R})$ with $\theta_t=u(t,L),$ $\intervala$ and
$$
\left\{
\begin{array}{ll}
u(t,\cdot)\in C^1(\mathbb{R}) \ \forall \ \interval. \\ \\
\int_{-n}^{n}\int_0^1 (\partial_xu(t,x))^2 dt dx <+\infty, \forall n \in \mathbb{N}  \ a.s.. \\ \\
\mathbb{E}\pqa \int_{\mathbb{R}}\pta  \int_0^1 (\partial_x u(t,x))^2 dt \ptc^2 dx+  \int_0^1 (u(t,0))^2 dt\pqc
 <+\infty. \\ \\
\mathbb{E}\pqa \int_0^1 (\partial_xu(t,L))^2 (D_tL)^2 dt
+\pta \int_0^1 (\partial_xu(t,L))^2dt \ptc \pta \int_0^1
(D_tL)^2 dt \ptc \pqc <+\infty.
\end{array}
\right.
$$
Suppose that $\cal{A}$ equals $\Theta(L).$ 
With the specifications above we have the following.
\begin{coro}\label{c4.31}
Let $b$ be a process in $L^2(\Omega \times [0,1]),$ such that $h(\cdot,L)+b$ belongs to the closure of $\cal A$ in $L^2(\Omega \times [0,1]).$ There exists an optimal process $\pi$ in $\cal{A}$ for the function 
$$
\theta \mapsto \mathbb{E}\pqa \int_0^1 \theta_td^-\pta W_t+\int_0^t b_sds\ptc-\frac{1}{2}\int_0^1 \theta_t^2dt \pqc
$$ if and only if $h(\cdot,L)+b$ belongs to $\cal{A}$ and $h(\cdot,L)+b=\pi$. 
\end{coro}
\begin{proof}
It is clear that $\cal{A}$ is a real linear set of measurable and with bounded paths processes verifying condition 1. of assumption \ref{a0}.  
Proposition 2.8 of \cite{LNN} shows that every $\theta$ in $\cal{A}$ is in $L^2(\Omega \times [0,1]),$ that it is  $W$-improperly forward integrable and that the improper integral belongs to $L^2(\Omega).$ In particular,  condition 2. of assumption \ref{a0} is verified. Furthermore,  the proof of theorem 3.2 of \cite{LNN} implicitly  shows that the process  
$
W-\int_0^\cdot h(t,L)dt,
$
is a $\mathcal{A}$-martingale. This implies that $W+\int_0^\cdot   b_t dt -\int_0^\cdot \gamma_t dt, $ with $\gamma=h(\cdot,L)+b,$ is an $\mathcal{A}$-martingale. The end of the proof follows then by proposition \ref{c3.18}.     
\end{proof}
\section{The market model}
We consider a market offering two investing possibilities in the time interval $[0,1].$ Prices of the two traded assets follow the evolution of two stochastic processes $\pta
S_t^0,\interval\ptc$ and $\pta S_t,\interval\ptc.$ We assume that
$$
S^0=\pta \exp(V_t),\interval\ptc,
$$
where $\pta V_t, 0\leq t\leq 1
\ptc$ is a 
positive process starting at zero with bounded variation,  
and $S$ is a continuous strictly positive process, with finite quadratic variation.

\begin{rema}
If  $V=\int_0^\cdot r_s ds,$ being $\pta r_t, \interval\ptc$ the short interest rate, $S^0$ represents the price process of the so called \textit{money market account}. 
Here we do not assume that $V$ is a riskless asset, being that assumption not necessary to develop our calculus. We only suppose that $S^0$ is \textit{less risky} then $S.$ 

Assuming that $S$ has a finite quadratic variation is not restrictive at least for  two reasons. 

Consider a market model involving an \textit{inside} trader: that means an investor having additional informations with respect to the \textit{honest} agent. Let $\mathbb{F}$ and $\mathbb{G}$ be the filtrations representing the information flow of the honest and the inside investor, respectively. Then it could be worthwhile to demand  the absence of  {\it{free lunches with vanishing risk}} (FLVR) among all simple $\mathbb{F}$-predictable strategies. Under the hypothesis of absence of (FLVR), by theorem 7.2, page 504 of \cite{DS}, $S$ is a semimartingale on the underlying probability space $(\Omega, P, \mathbb{F}).$  On the other hand $S$ could fail to be a $\mathbb{G}$-semimartingale, since (FLVR) possibly exist for the insider. Nevertheless, the inside investor is still allowed to suppose that $S$ has finite quadratic variation thanks to proposition \ref{p3.4}.

Secondly, as already specified in the introduction, if we want to include $S$  as a \textit{self-financing}-portfolio, we have to require that $\int_0^\cdot Sd^-S$ exists. This is equivalent to assume that $S$ has finite quadratic variation, see proposition 4.1 of \cite{RV2}.
\end{rema}
\subsection{Portfolio strategies}
We assume the point of view of an investor whose flow of information
is modeled by a filtration \filtrationg$=\pta
\mathcal{G}_t\ptc_{t\in[0,1]}$ of $\cal F,$ which satisfies the usual assumptions.

We denote with $C^-_{b}([0,1))$ the set of processes  which have paths being left continuous and bounded on each compact set of $[0,1).$

\begin{defi}\label{d4.4}
A \textbf{portfolio strategy} is a couple of $\mathbb{G}$-adapted processes $\phi=\pta\pta h^0_t,h_t\ptc, 0\leq t<1\ptc.$ 
The market value $X$ of the portfolio strategy $\phi$ is the so called \textbf{wealth process}
$
X=h^0S^0+hS.
$
\end{defi}

We stress that there is no point in defining the portfolio strategy at the end of the trading period, that is for $t=1.$  Indeed, at time $1,$ the  agent has to liquidate his portfolio.

\begin{defi}\label{d5.3}
A portfolio strategy $\phi=\pta h^0,h\ptc$ 
is \textbf{self-financing} if  both $h^0$ and $h$ belong to $C^-_{b}([0,1)),$ the process $h$ is locally $S$-forward integrable and its wealth process $X$ verifies 
\beqn \label{sfc}
X=X_0+\int_0^\cdot  h^0_tdS^0_t + \int_0^\cdot h_t d^-S_t.
\eeqn
\end{defi}

The interpretation of the first two items is straightforward: $h^0$ and $h$ represent, respectively, the number of shares of $S^0$ and $S$ held in the portfolio; $X$ is its market value. The self-financing condition (\ref{sfc}) seems to be an appropriate formalization of the intuitive idea of trading strategy not involving exogenous sources of money. Among its justifications we can include the following ones.

As already explained in the introduction, the discrete time version of condition  (\ref{sfc}) reads as the classical self-financing condition. 
Furthermore, if
    $S$ is a $\mathbb{G}$-semimartingale, forward integrals of $\mathbb{G}$-adapted processes with left continuous and bounded
 paths, agree with  classical It\^o integrals, see proposition \ref{p3.6} and \ref{p3.4}.

In the sequel we will choose as \textit{num\'eraire} the positive process $S^0.$ That means that prices will be expressed in terms of $S^0.$  We will denote with $\widetilde{Y}$
the  value of a stochastic process $(Y_t,\interval)$ \textit{discounted} with respect to $S^0:$
$
\widetilde{Y_t}={Y_t}({S^0_t})^{-1},$ for every $\interval.$

The following lemma shows that, as well as in a semimartingale model, a portfolio strategy which is self-financing is uniquely determined by its initial value and the process representing the number of shares of $S$ held in the portfolio.
We remark that previous definitions and considerations can be made without supposing that the investor is able to observe prices of $S$ and $S^0.$ 
However, we need to make this hypothesis for the following characterization of self-financing portfolio strategies.    
\begin{assu}\label{a5.4}
>From now on we suppose that $S$ and $S^0$ are $\mathbb{G}$-adapted processes.
 \end{assu}
\begin{lemm}\label{l5.5}
Let $\phi=\pta\pta h^0_t,h_t\ptc, \intervala\ptc$ be a couple of $\mathbb{G}$-adapted processes in $C^-_{b}([0,1)).$ Suppose that $h$ is locally $S$-forward integrable.  Then the portfolio strategy $\phi$ is self-financing
 if and only if its discounted wealth process $\widetilde{X}$ verifies
\beqn \label{wealth0}
\widetilde{X}=X_0-\int_0^\cdot e^{-V_t}h_tS_tdV_t+ \int_0^\cdot e^{-V_t}d^- \int_0^t h_sd^-S_s.
\eeqn
On the other hand, let $\pta h_t,\intervala\ptc$ be a $\mathbb{G}$-adapted  process in $C^-_{b}([0,1)),$ which is locally $S$-forward integrable, and  $X_0$ be a  $\mathcal{G}_0$-random variable. Then
the couple $$\phi=\pta \pta (X_t-h_tS_t)(S^0_t)^{-1},h_t\ptc,\intervala \ptc,$$ with $X$ defined as in equality (\ref{wealth0}), is a self-financing portfolio strategy with wealth process $X.$
\end{lemm}
\begin{proof}
Regarding the first part of the statement we observe that corollary \ref{c2.15} and equality $X=h^0B+h S$ 
imply the equivalence between equalities (\ref{wealth0}) and (\ref{sfc}).

Let $h,$ $X_0$ and $X$ be as in the second part of the statement. It is clear that $h^0=\pta \pta X_t-h_tS_t\ptc (S^0_t)^{-1},\intervala\ptc$ is $\mathbb{G}$-adapted and belongs to $C^-_{b}([0,1)).$ By construction, the wealth process corresponding to the strategy $\phi=(h^0,h)$ is equal to $X.$ The conclusion follows by the first part of the statement.
\end{proof}

\begin{rema}\label{r5.6}
Suppose that $(h^0,h)$ is a self-financing portfolio strategy with $h$ locally forward integrable with respect to $\widetilde{S}.$ Corollary \ref{c2.15} and previous lemma imply that its discounted wealth process $\widetilde{X}$ can be also be rewritten in the following way
$$
\widetilde{X}=X_0+\int_{0}^\cdot h_td^-\widetilde{S}_t+R,
$$
with 
$$
R=\int_0^\cdot e^{-V_t}d^- \int_0^t h_sd^-S_s-\int_{0}^\cdot h_td^-\int_0^t e^{-V_s}d^- S_s. 
$$
\end{rema}

Lemma \ref{l5.5} leads to conceive the following definition.
\begin{defi} \label{d5.4}
\begin{enumerate}
\item
A \textbf{self-financing portfolio} is a couple $\pta X_0, h\ptc$ of a
$\mathcal{G}_0$-measurable random variable $X_0$, and  a process $h$ in $C_b^-([0,1))$ which is $\mathbb{G}$-adapted and locally $S$-forward integrable. 
\item
The discounted wealth process $\widetilde{X}$ of the self financing portfolio $\pta X_0, h\ptc,$ and the number of shares $h^0$ of $S^0$ held in that portfolio are  given by
$$
\left\{
\begin{array}{ll}
\widetilde{X}=X_0-\int_0^\cdot e^{-V_t}h_tS_tdV_t+ \int_0^\cdot e^{-V_t}d^- \int_0^t h_sd^-S_s\\
h^0=(X-hS)(S^0)^{-1}.
\end{array}
\right.
$$
\item In the sequel we let us employ the term \textit{\textbf{portfolio}} to denote the process $h,$ in a self-financing portfolio, representing the number of shares of $S$ held. Without further specifications the initial wealth of an investor will be assumed to be equal to zero. 
\end{enumerate}
\end{defi}

Lemma \ref{l5.5} and remark \ref{r5.6} immediately  imply the following.
   
\begin{coro} \label{l5.6} Let $(X_0,h)$ be a self-financing portfolio such that $h$ is locally $\widetilde{S}$-forward integrable  and  $\int_0^\cdot e^{-V_t}d^-\int_0^t h_s d^-S_s=\int_0^\cdot h_t d^-\int_0^t e^{-V_s}d^-S_s.$ Then its discounted value $\widetilde{X}$ verifies the equality
$
\widetilde{X}=X_0+\int_0^\cdot hd^-\widetilde{S}.
$ 
\end{coro}

\begin{rema}
If $S$ is a $\mathbb{G}$-semimartingale, the hypothesis required on $h$ in previous remark is always verified. Indeed, forward integrals coincide with classical It\^o integrals for which the associative property holds true, see proposition \ref{p3.4}. 

Some conditions to insure the existence of \textit{chain-rule} formulae, when the semimartingale property of the integrator process fails to hold, can be found in section \ref{s3}. For more informations about this topic we also refer to  \cite{FR} and \cite{ER}.
\end{rema}

\begin{assu}\label{a5.9} 
We assume the existence of a real linear space of portfolios ${\bf{\cal A}}$, that is of $\mathbb{G}$-adapted processes $h$ belonging to $C^-_b([0,1)),$ which are locally $S$-forward integrable. The set $\cal{A}$ will represent the set of all \textbf{admissible strategies} for the investor.
\end{assu}

We proceed furnishing examples of sets behaving as the set $\cal{A}$ in assumption \ref{a5.9}. They correspond to the examples discussed in section \ref{s3}.

\subsubsection{Admissible strategies via It\^o fields} \label{s5.1.1}
We refer the reader to subsection \ref{s4.0} for notations and definitions.

The following proposition is a straightforward consequence of proposition \ref{3.8}.

\begin{prop}\label{p5.11}
Let $\cal{A}$ be the set of processes $(h_t, \intervala)$  such that for every $\intervala$ the process in $hI_{[0,t]}$ belongs to $\mathcal{S}(\mathcal{C}^1_S(\mathbb{G})).$ Then $\cal{A}$ is a real linear space satisfying the hypotheses of assumption \ref{a5.9}.
\end{prop}

\subsubsection{Admissible strategies via Malliavin calculus}\label{5.25}
For this example we refer to subsection \ref{s4.1}. We recall that there,  $W$ was a real valued Wiener process defined on the canonical probability space $\pta \Omega, \mathbb{F}, \mathcal{F},P\ptc.$ 
Regarding the price of $S$ we make the following assumption.
\begin{assu}
We suppose that
$
S=S_0\exp \pta \int_0^\cdot \sigma_t dW_t+\int_0^\cdot\pta \mu_t-\frac{1}{2}\sigma_t^2\ptc dt\ptc,$ 
where $\mu$ and $\sigma$ are $\mathbb{F}$-adapted, $\mu$ belongs to $L^{1,q}$ for some $q>4,$ $\sigma$ has bounded and left continuous paths, it  belongs 
$L^{1,2}_-\cap L^{2,2}$ and the random variable 
\beqn \nonu
\sup_{t\in [0,1]} \pta\vaa \sigma_t \vac+\sup_{s\in[0,1]} \vaa D_s \sigma_t\vac 
\sup_{s,u\in[0,1]} \vaa D_s D_u\sigma_t\vac \ptc \eeqn 
is bounded.
\end{assu}

\begin{rema}\label{r5.15}
By remark of page 32, section 1.2 of \cite{N} 
 $\sigma$ is in $L^{1,2}_-$ and $D^-\sigma=0$.
\end{rema}

Using remarks \ref{r4.1} and \ref{r4.3}, lemma \ref{l4.4} and lemma \ref{l4.6}, it is not difficult to prove that
the process 
$
\log\pta S\ptc
$ 
belongs to $L^{1,q}_-.$

\begin{prop} \label{p5.16} Let $\mathcal{A}$ be the set of all $\mathbb{G}$-adapted processes $h$ in $C^-_b([0,1)),$ such that for every $0\leq t<1,$ the process $hI_{[0,t]}$ belongs to $L^{1,p}_{-},$ for some $p>4.$  
Then $\cal{A}$ is a real linear space satisfying the hypotheses of assumption \ref{a5.9}.
\end{prop}

\begin{proof}
Let $h$ be in $\cal{A}.$ We set 
$
A=\log(S)-\log(S_0)+\frac{1}{2}\int_0^\cdot\sigma_t^2dt=\int_0^\cdot \sigma_tdW_t+\int_0^\cdot \mu_tdt.
$
We recall that, thanks to lemma \ref{l2.1}, for every $0\leq t <1,$ $hI_{[0,t]}$ is $S$-forward integrable if and only if $hI_{[0,t]}S$ is forward integrable with respect to $A.$ Let $0\le t<1,$ be fixed.
Each component of the vector process $u=\pta hI_{[0,t]},\log(S)\ptc$ belongs to $L^{1,p}_-$ for some $p>4$ and it has left continuous and bounded paths. We can thus apply proposition \ref{p4.15} to state that $hI_{[0,t]}S$ is forward integrable with respect to $\int_0^\cdot \sigma_t dW_t.$ 
This implies that $hI_{[0,t]}S$ is $A$-forward integrable. Letting $t$ vary in $[0,1)$
we find that $h$ is $S$-improperly integrable and we get the end of the proof. \end{proof}

\subsubsection{Admissible strategies via substitution} \label{e5.21}
We consider the setting of subsection \ref{s4.2}. More precisely, we assume the existence of a filtration $\mathbb{F}=\pta \mathcal{F}_t\ptc_{t\in[0,1]}$ on $\pta \Omega, \mathcal{F}, P\ptc,$ with $\mathcal{F}_1=\mathcal{F},$ and of an $\mathcal{F}$-measurable random variable $L$ with values in $\mathbb{R}^d$ such that
$\mathcal{G}_t=\pta \mathcal{F}_{t} \vee \sigma(L)\ptc
$, for every $\interval.$ We suppose that $\mathbb{G}$ is right continuous. We assume that $S$ and $S^0$ are $\mathbb{F}$-adapted, and that $S$ is an   $\mathbb{F}$-semimartingale.

We observe that this situation arises when the investor trades as an \textit{\textbf{insider}}, that is having an extra information about prices, at time $0,$ represented by the random variable $L.$

\begin{prop}
Let $\cal{A}$ be the set of processes $h$ such that, for every $0\leq t <1,$ the process $hI_{[0,t]}$ belongs to $\mathcal{S}(\mathcal{A}^{p,\gamma}(L)),$ for some $p>1$ and $\gamma>d.$ Then $\cal{A}$ satisfies the hypotheses of  assumption \ref{a5.9}.  
\end{prop}

\begin{proof}
Processes in $\cal{A}$ are clearly $\mathbb{G}$-adapted and in $C^-_b([0,1)).$ 
The end of the proof is a consequence of proposition \ref{p3.35} and remark \ref{r5.2}.
\end{proof}

The following lemma shows that is not so reductive to restrict the class of possible portfolio strategies to the collection of sets $\pta \mathcal{S}\pta \mathcal{A}^{p,\gamma}(L)\ptc, p>1,\gamma>d \ptc$.

\begin{lemm} \label{l4.13}
Let $\pta\pi_t, 0\leq t < 1 \ptc$ be a bounded $\mathcal{P}^{\mathbb{G}}$-measurable process. Then there exists a $\mathcal{P}^{\mathbb{F}}\otimes \mathcal{B}(\mathbb{R})$-measurable function $\pta h(t,x),\intervala, x \in \mathbb{R}^{d}\ptc $,  such that $\pi=h(\cdot,L),$ almost surely.
\end{lemm}

\begin{proof}
Define $L^{0,\mathcal{P}^{\mathbb{F}}}$ as the set of all functions $(h(t,x), \intervala, x\in \mathbb{R}^d)$ which are $\mathcal{P}^{\mathbb{F}} \otimes \mathcal{B}(\mathbb{R}^{d})$-measurable.
Set
$$
\mathcal{M}=\pga u: \Omega\times [0,1)\rightarrow \mathbb{R}, \left|\right.  \exists h \in L^{0,\mathcal{P}^{\mathbb{F}}}, \mbox{ s.t. }  h(\cdot,L)=u \mbox{ a.s. } \pgc.
$$
The set $\mathcal{M}$ is a \textit{monotone vector space}, see definition in chapter 1 of \cite{P}. Indeed, it is a linear vector space of bounded real functions containing all constants and,
if $\pta u_n\ptc_{n\in \mathbb{N}}$  is an increasing sequence of positive random elements in $\mathcal{M},$ with $u=\sup_{n\in \mathbb{N}}u_n$ bounded, then $u$ belongs to $\mathcal{M}.$ In fact $
h=sup_{n}h_n
$
is still in $L^{0,\mathcal{P}^{\mathbb{F}}}$ and $u=h(\cdot,L).$
Consider the set $\mathcal{S}^{\mathbb{G}}$ of all $\mathcal{P}^{\mathbb{G}}$-measurable processes of the form
$
u=I_{\pga 0\pgc }h_0(L)f_0+\sum_{i=0}^{k-2} I_{(t_{i},t_{i+1}]}h_i(L)f_i+I_{(t_{k-1},1)}h_{k-1}(L)f_{k-1},
$
where $0=t_0<t_1<...<t_k=1$, and $h_i$ is $\mathcal{B}(\mathbb{R})$-measurable and bounded, $f_i$ is $\mathcal{F}_{t_i}$-measurable and bounded, for every $i=0,...,k.$ It is clear that  $\mathcal{S}^{\mathbb{G}}$ is stable with respect to multiplication. Moreover  $\sigma\pta\mathcal{S}^{\mathbb{G}}\ptc$ contains the $\sigma$-algebra generated by all bounded and $\mathcal{P}^{\mathbb{G}}\otimes \mathcal{B}(\mathbb{R}^d)$-measurable function. We can thus apply theorem $8$ of \cite{P} to get the result.
\end{proof}

\subsection{Completeness and arbitrage: $\mathcal{A}$-martingale measures}

\begin{defi}
Let $h$  be a self financing portfolio in $\cal{A},$ which is $S$-improperly forward integrable and $X$ its wealth process. Then $h$ is an \textbf{$\mathcal{A}$-arbitrage} if  $X_1=\lim_{t\rightarrow 1}X_t$ exists almost surely, $P(\pga X_1 \geq 0\pgc )=1$ and  $P(\pga X_1>0\pgc)>0.$ 
\end{defi} 

\begin{defi}
If there are no $\mathcal{A}$-arbitrages we say that the market is \textbf{$\mathcal{A}$-arbitrage free}.
\end{defi}
\begin{defi}
A probability measure $Q \sim  P$ is said \textbf{$\mathcal{A}$-martingale measure} if under $Q$ the process $\widetilde{S}$ is an $\mathcal{A}$-martingale according to definition \ref{d3.1}.
\end{defi}

We will need the following assumption.
\begin{assu} \label{a}
Suppose that for all $h$ in $\cal{A}$  the following conditions hold.
\begin{enumerate}
\item The process $e^{V}h$ belongs to $\mathcal{A}.$
\item $h$ is $\widetilde{S}$-improperly forward integrable and 
\beqn \label{e6}
\int_0^\cdot e^{-V_t}d^-\int_0^t h_sd^-S_s=\int_0^\cdot h_t e^{-V_t}d^-S_t=\int_0^\cdot h_td^-\int_0^t e^{-V_s}d^-S_s.
\eeqn
\end{enumerate}
\end{assu}

For the following proposition the reader should keep in mind the notation in equality (\ref{AX}). We omit its proof which is a direct application of corollary \ref{c4.16}.

\begin{prop} \label{p5.27}
Let  $\mathcal{A}=\mathcal{A}_S.$  Suppose that $d\pqa S\pqc_t=\sigma(t,S_t)^2dt,$ where $\sigma$ satisfies assumption \ref{a4.15}. If there exists a $\cal A$-martingale measure then  the law of $\widetilde{S}_t$ is absolutely continuous with respect to Lebesgue measure, for every $\interval.$ 
\end{prop}

\begin{prop} \label{p5.28}
Under assumption \ref{a}, if there exists an $\mathcal{A}$-martingale measure $Q,$ the market is $\mathcal{A}$-arbitrage free.
\end{prop}
\begin{proof}
Suppose that $h$ is an $\cal{A}$-arbitrage. Since $\widetilde{S}$ is an $\mathcal{A}$-martingale under $Q,$ using corollary \ref{l5.6} we find $\mathbb{E}^Q[\widetilde{X}_1]=\mathbb{E}^Q[\int_0^1 h_td^-\widetilde{S}_t]=0$. This contradicts the arbitrage condition $Q(\pga X_1>0\pgc)>0.$
\end{proof}

The proposition which follows characterizes the set of all $\mathcal{A}$-martingale measures.

\begin{prop}
Under assumption \ref{a} the process $\widetilde{S}$ is an $\mathcal{A}$-martingale under $Q,$ if and only if 
the process $
S-\int_0^\cdot S_t dV_t
$
is an $\mathcal{A}$-martingale under $Q$.
\end{prop}
\begin{proof}
If $h$ is in $\mathcal{A}$ by assumption \ref{a} we have 
\beqn \nonu
\mathbb{E}^Q\pqa \int_0^\cdot h_t d^- \pta  S_t-\int_0^t S_s dV_s \ptc\pqc  &=&\mathbb{E}^Q\pqa\int_0^\cdot (h_te^{V_t})e^{-V_t} d^- \pta  S-\int_0^t S_s dV_s \ptc \pqc\\ \nonu
&=& \mathbb{E}^Q\pqa\int_0^\cdot (h_te^{V_t})d^- \int_0^t e^{-V_s}d^-S_s \pqc\\ \nonu &+&\mathbb{E}^Q\pqa \int_0^\cdot (h_te^{V_t})d^-\int_0^t S_s de^{-V_s}  \pqc  \\ \nonu &=&\mathbb{E}^Q\pqa\int_0^\cdot (h_te^{V_t})d^-\widetilde{S_t}\pqc=0.
\eeqn
\end{proof}

We proceed discussing \textit{completeness}.
\begin{defi}
A  \textbf{contingent claim} is an $\cal{F}$-measurable random variable.
 $\mathcal{L}$ will be a set of $\cal{F}$-measurable random variables; it will represent all the {contingent claims} the investor is interested in.
\end{defi}

\begin{defi} \begin{enumerate}
\item A contingent claim $C$ is said \textbf{$\mathcal{A}$-attainable} if there exists a self financing portfolio $(X_0,h)$ with $h$
in $\mathcal{A},$ which is $S$-improperly forward integrable,  
such that the corresponding wealth process $X$ verifies
$\lim_{t\rightarrow 1}X_t=C,$ almost surely. The portfolio $h$ is said the \textbf{replicating} or \textbf{hedging} portfolio for $C,$ $X_0$ is said the \textbf{replication price} for $C.$
\item The market is said to be \textbf{$\pta \mathcal{A},\mathcal{L}\ptc$-complete} if every contingent claim in $\cal{L}$ is attainable trough a portfolio in $\mathcal{A}.$
\end{enumerate}
\end{defi}

\begin{assu}\label{A1}
For every $\mathcal{G}_0$-measurable random variable $\eta,$ and $h$ in $\cal{A}$ the process $u=h\eta,$ belongs to $\cal{A}.$ 
\end{assu}

\begin{prop}\label{p5.33}
Suppose that the market is $\cal{A}$-arbitrage free, and that assumption $\ref{A1}$ is realized. Then the replication price of an attainable contingent claim is unique. 
\end{prop}
\begin{proof}
Let $(X_0,h)$ and $(Y_0,k)$ be two replicating portfolios for a contingent claim $C,$ with $h$ and $k$ in $\cal{A},$ and wealth processes $X$ and $Y$, respectively. We have to prove that $$
P\pta \pga X_0-Y_0\neq 0\pgc \ptc=0.
$$ 
Suppose, for instance, that $P\pta X_0-Y_0>0\ptc \neq 0.$ We set $A=\pga  X_0-Y_0>0\pgc.$ By assumption \ref{A1}, $I_{A}(k-h)$ is a portfolio in $\cal{A}$ with wealth process  $I_A (Y_t-X_t).$ Since both $(X_0,h)$ and $(Y_0,k)$ replicate $C,$ $\lim_{t\rightarrow 1} I_{A}(Y_t-X_t)=I_{A}(X_0-Y_0),$ with $P(\pga I_A(X_0-Y_0>0)\pgc )>0.$ Then  $I_A(k-h)$ is an $\cal{A}$-arbitrage and this contradicts the hypotheses.   
\end{proof}

\begin{prop}\label{p5.30} Suppose that there exists an $\mathcal{A}$-martingale measure $Q$.  Then the following statements are true.
\begin{enumerate}
\item Under assumptions \ref{a} and  \ref{A1}, the replication price of an $\cal{A}$-attainable contingent claim $C$ is unique and equal to $\mathbb{E}^Q\pqa \widetilde{C}\left|\right.\mathcal{G}_0\pqc.$
\item Let $\mathcal{G}_0$ be trivial. If  $Q$ and $Q_1$ are two $\mathcal{A}$-martingale measures, then $\mathbb{E}^Q[ \widetilde{C}]=\mathbb{E}^{Q_1}[ \widetilde{C}],$ for every $\cal{A}$-attainable contingent claim $C$. In particular, if  the market is $(\cal{A},\cal{L})$-complete and $\mathcal{L}$ is an algebra, all $\mathcal{A}$-martingale measures coincide on the $\sigma$-algebra generated by all bounded discounted contingent claims in $\mathcal{L}.$ 
\end{enumerate}
\end{prop}
\begin{proof}
Let $(X_0,h)$ be a replicating $\cal{A}$-portfolio for $C$.
By corollary \ref{l5.6} 
$$
\mathbb{E}^{Q}\pqa \widetilde{C}\left|\right. \mathcal{G}_0\pqc=X_0+\mathbb{E}^{Q}\pqa \int_0^1 h_td^-\widetilde{S}_t \left|\right.\mathcal{G}_0\pqc.
$$
We observe that 
$
\mathbb{E}^{Q}\pqa \int_0^1 h_td^-\widetilde{S}_t \left|\right.\mathcal{G}_0\pqc=0.
$
In fact, if $\eta$ is  a $\mathcal{G}_0$-measurable random variable, then, thanks to  assumption  \ref{A1}, $\eta h$ belongs to $\mathcal{A},$ so as to have
$\mathbb{E}^Q\pqa \pta \int_0^1 h_td^-\widetilde{S}_t\ptc \eta \pqc=\mathbb{E}^Q\pqa  \int_0^1 \eta h_td^-\widetilde{S}_t\pqc=0.$  This implies point 1.

If $\mathcal{G}_0$ is trivial, we deduce that, if $Q$ and $Q_1$ are two $\mathcal{A}$-martingale measures,  $\mathbb{E}^{Q}[ \widetilde{C}]=\mathbb{E}^{Q_1}[\widetilde{C}],$ for every $\mathcal{A}$-attainable contingent claim. The proof of the last point is then an application of theorem 8, chapter 1 of \cite{P}. 
\end{proof}
\subsection{Hedging}
In this part of the paper we price contingent claims via partial differential equations. In particular we show robustness of Black-Scholes formula for \textit{European} and \textit{Asian} contingent claims within a non-semimartingale model. 

The following proposition generalizes a result obtained in a slight different form in \cite{Z}, when the process $S$ is supposed to be the sum of a Wiener process and a continuous process with zero quadratic variation.

We suppose here that $d\pqa S \pqc_t=\sigma(t,S_t)^2 S_t^2dt$ and $dV_t=rdt,$ with $ r>0$ and $\sigma:[0,1]\times (0,+\infty)\rightarrow \mathbb{R}.$   

\begin{prop} \label{p5.31}
Let $\psi$ be a function in $C^0(\mathbb{R})$. Suppose that there exists $\pta v(t,x), \interval, x\in \mathbb{R}\ptc$ of class $C^{1,2}([0,1)\times \mathbb{R})\cap C^0([0,1]\times \mathbb{R}),$ which is a  solution of the following Cauchy problem
\beqn \label{e14}
\left\{
\begin{array}{lll}
\partial_tv(t,y)+\frac{1}{2}(\widetilde{\sigma}(t,y))^2y^2\partial_{yy}^{(2)}v(t,y)&=&0 \quad \mbox{ on }  [0,1)\times{\mathbb{R}}
\\
v(1,y)&=&\widetilde{\psi}(y),
\end{array}
\right.
\eeqn
where 
$$
\left\{
\begin{array}{ll}
\widetilde{\sigma}(t,y)=\sigma(t,ye^{rt})&\quad \forall (t,y) \in [0,1]\times\mathbb{R},
\\ \widetilde{\psi}(y)=\psi(ye^r)e^{-r}&\quad \forall y \in \mathbb{R}.  
\end{array}
\right. 
$$
Set 
$$
h_t=\partial_yv(t,\widetilde{S}_t),\quad \intervala,  \quad X_0=v(0,S_0).
$$
Then $(X_0,h)$ is a self-financing portfolio replicating the contingent claim $\psi(S_1).$
\end{prop}
\begin{proof}
Assumption \ref{a5.4} tells us that $h$ is a $\mathbb{G}$-adapted process in $C^-_b([0,1)).$  By proposition \ref{p2.1}, $h$ is locally $\widetilde{S}$-forward integrable. Combining lemma \ref{l2.1} and proposition \ref{p2.1}, it is possible to prove that  
$$
\int_0^\cdot e^{-V_t}d^-\int_0^t h_s d^-S_s=\int_0^\cdot h_t e^{-V_t}d^-S_t=\int_0^\cdot h_t d^-\int_0^t e^{-V_s}d^-S_s.
$$ 
Similar arguments were used in \cite{FR},  corollary 23.
Corollary \ref{l5.6} implies then that its discounted wealth process  verifies 
\beqn \label{e19}
\widetilde{X}=X_0+\int_0^\cdot hd^-\widetilde{S}.
\eeqn
On the other hand by point 2. of proposition \ref{3.8} 
\beqn \label{e20}
[ \tilde{S}]
=\int_0^\cdot \widetilde{S}^2_s \widetilde{\sigma}(s,\widetilde{S}_s)^2ds.
\eeqn
Applying  proposition \ref{p2.1}, recalling equation (\ref{e14}), equalities (\ref{e19}) and (\ref{e20}) we find that  
$$
\widetilde{X}_t=v(t,\widetilde{S}_t), \quad \forall \intervala.
$$
In particular $X_0+\lim_{t\rightarrow 1} \int_0^t h_sd^-\widetilde{S}_s$ exists finite and coincides with $v(1,\widetilde{S}_1)=\widetilde{\psi}(\widetilde{S_1})=\psi(S_1)e^{-r}.$
\end{proof}

\begin{rema}
In particular, under some minimal regularity assumptions on  $\sigma$ and  no degeneracy, the market  is $(\mathcal{A}_S,\mathcal{L})$-complete, if $\mathcal{L}$ equals the set of all contingent claims of type $\psi(S_1)$ with $\psi$ in $C^0(\mathbb{R})$ with linear growth.
\end{rema}

The result of proposition \ref{p5.31} can also be adapted to hedge Asian contingent claims, that is contingent claims  $C$ depending on the mean of $S$ over the traded period: $C=\psi\pta  \frac{1}{S_1}\pta \int_0^1 S_t dt\ptc\ptc,$ for some $\psi$ in $C^0(\mathbb{R}).$  
\begin{prop}
Suppose that $\sigma(t,x)=\sigma,$ for every  $(t,x)$ in $[0,1]\times \mathbb{R},$ for some $\sigma>0.$ Let $\psi$ be a function in $C^0(\mathbb{R})$  and $v(t,y)$ a continuous solution of class $C^{1,2}([0,1)\times  \mathbb{R})\cap C^{0}([0,1]\times \mathbb{R})$ of the following Cauchy problem
\beqn \nonu
\left\{
\begin{array}{lll}
\frac{1}{2}\sigma^2 y^2 \partial_{yy}^{(2)}v(t,y)+ (1-ry)\partial_{y}v(t,y)+\partial_tv(t,y)&=&0, \quad \mbox{ on }  [0,1)\times{\mathbb{R}}
\\
v(1,y)&=&\psi(y).
\end{array}
\right.
\eeqn
Set  $Z_t=\int_0^t S_sds-K,$ 
for some $K>0,$ $X_0=v(0,\frac{K}{S_0})S_0$ and $h_t=v(t,\frac{Z_t}{S_t})-\partial_y v(t,\frac{Z_t}{S_t})\frac{Z_t}{S_t},$ for all $0\leq t \leq 1.$ Then $\pta X_0,h \ptc$ is a self-financing portfolio which replicates the contingent claim $\psi \pta \frac{1}{S_1}\pta \int_0^1 S_tdt-K\ptc\ptc S_1.$
\end{prop}

\begin{proof}
We set $\xi_t=\frac{Z_t}{S_t}, \interval$. Applying proposition \ref{p2.1} to the function $u(t,z,s)=
v(t,\frac{z}{s}e^{-rt})s$ and using the equation fulfilled by $v$ we can 
 expand the process  $(e^{-rt}v(t,\xi_t)S_t, \intervala)$ as  follows:
\beqn \label{f19}
u(t,Z_t,\widetilde{S}_t)=v\pta t,\xi_t\ptc\tilde{S}_t=v\pta 0,\xi_0\ptc S_0+\int_0^t h_td^-\widetilde{S}_t.
\eeqn 
By arguments which are similar to those used in the proof of previous proposition, it is possible to show that $h$ is a self-financing portfolio and that (\ref{f19}) implies that $u(t,Z_t,\widetilde{S}_t)=\widetilde{X}_t$ for every $\intervala.$ Therefore $\lim_{t \rightarrow 1}\widetilde{X}_t$ is finite and equal to $\psi\pta \xi_1\ptc S_1 e^{-r}.$ This concludes the proof. 
\end{proof}

\subsection{Utility maximization}
\subsubsection{Formulation of the problem}

We consider the problem of maximization of expected utility from terminal wealth starting
from initial capital $X_0> 0,$ being $X_0$ a $\mathcal{G}_0$-measurable random variable. We define the function $U(x)$ modeling the
utility of an agent with wealth $x$ at the end of the trading period. The function $U$ is supposed to be of class $C^2((0,+\infty)),$ strictly increasing, with $U'(x)x$ bounded.

We will need the following assumption.

\begin{assu} \label{a5.41}
The utility function $U$ verifies  
$
\frac{U^{''}(x)x}{U'(x)}\leq -1, \quad \forall x>0.
$ 
\end{assu}

A typical example of function $U$ verifying assumption \ref{a5.41} is $U(x)=\log(x).$

We will focus on portfolios with strictly positive value. As a consequence of this, before starting analyzing the problem of maximization, we show how it is possible to construct portfolio strategies when only positive wealth is allowed.

\begin{defi} \label{d5.7} For simplicity of calculation we introduce the process 
$$
A=\log(S)-\log(S_0)+\frac{1}{2}\int_0^\cdot\frac{1}{S^2_t}d\pqa S\pqc_t.
$$
\end{defi}

\begin{lemm} \label{l4.3}
Let $\theta=\pta \theta_t, \intervala \ptc$ be a $\mathbb{G}$-adapted process in $C^-_{b}([0,1))$ such that 
\begin{enumerate}
\item $\theta$ is $A$-improperly forward integrable. 
\item The process $A^\theta=\int_0^\cdot \theta_sd^-A_s$ has finite quadratic variation.
\item If $X^\theta$ is the process defined by  
$$
X^\theta=X_0\exp \pta \int_0^\cdot \theta_t d^-A_t+\int_0^\cdot  \pta 1-\theta_t\ptc dV_t-\frac{1}{2}\pqa A^\theta \pqc \ptc,
$$
then
$\int_0^\cdot X^\theta_t \theta_td^-A_t$ and $\int_0^\cdot X^\theta_t d^-\int_0^t \theta_s d^-A_s
$ improperly exist  and 
\beqn \label{c10} 
\int_0^\cdot X^\theta_t \theta_td^-A_t=\int_0^\cdot X^\theta_t d^-\int_0^t \theta_s d^-A_s.
\eeqn
\end{enumerate}
Then the couple $\pta X_0, h\ptc$, with  $h_t=\frac{\theta X^\theta_t}{S_t},$ $\intervala,$
is a self-financing portfolio with strictly positive wealth $X^\theta.$ In particular,  $\lim_{t\rightarrow 1} X^\theta_t=X^\theta_1$ exists and it is strictly positive.  
\end{lemm}

\begin{proof}
Thanks to lemma \ref{l2.1} $h$ is locally $S$-forward integrable and 
$
\int_0^\cdot h_t d^-S_t=\int_0^\cdot \theta_t X^\theta_td^-A_t.
$
Applying corollary \ref{c2.15}, proposition \ref{p2.1}, and using hypothesis  3., $\widetilde{X^\theta}$ can be rewritten in the following way:
\beqn \nonu
\widetilde{X^\theta_t}&=&X_0+\int_0^\cdot \widetilde{X^\theta_t}d^-\int_0^t \theta_s d^-A_s-\int_0^\cdot \widetilde{X^\theta_t}\theta_t dV_t\\ \label{e22}
&=&
X_0+\int_0^\cdot e^{-V_t}d^- \int_0^t h_sd^-S_s-\int_0^\cdot e^{-V_t}h_tS_tdV_t.
\eeqn
Remark  \ref{r5.6} tells us that $X^\theta$ is the wealth of the self-financing  portfolio $\pta X_0, h\ptc.$
\end{proof}

\begin{rema}
The process $\theta$ in previous lemma represents the \textit{\textbf{proportion}} of wealth invested in $S.$
\end{rema}
\begin{rema}
Let $\theta$ be as in lemma \ref{l4.3}. 
Then, for every $0\leq t < 1,$ $X$ is, indeed, the unique solution, on $[0,t],$ of equation
\beqn \nonu
X^\theta=X_0+\int_0^\cdot X^\theta_td^-\pta \int_0^t \theta_sd^-A_s +\int_0^t (1-\theta_s)dV_s-\frac{1}{2}\pqa A^\theta \pqc_t\ptc.
\eeqn
In fact, uniqueness is insured by corollary 5.5 of \cite{RV2}.
It is important to highlight that, without the assumption on $\theta$ regarding the  chain rule in equality (\ref{c10}),  we cannot conclude that $X^\theta$ solves equation (\ref{e22}).
However we need to require that $X^\theta$ solves the latter equation to interpret it  as  the value of a portfolio whose proportion invested in $S$ is constituted by $\theta.$ 
In the sequel we will construct, in some specific settings, classes of processes defining proportions of wealth as in lemma \ref{l4.3}. We will consider, in particular, two cases already contemplated in \cite{BO} and \cite{LNN}. Our definitions of those sets will result more complicated than the ones defined in the above cited papers. This happens because, in those works, the chain rule problem arising when the forward integral replaces the classical It\^o integral is not clarified.
\end{rema}

\begin{assu} \label{a5.40}
We assume the existence of a real linear space $\cal{A}^+$ of $\mathbb{G}$-adapted processes $(\theta_t,\intervala)$ in $C_b^-([0,1)),$  such that 
\begin{enumerate}
\item $\theta$ verifies condition 1., 2. and 3. of lemma \ref{l4.3}, and $\pqa A^\theta\pqc=\int_0^\cdot \theta_t^2d\pqa A\pqc_t.$
\item $\theta I_{[0,t]}$ belongs to $\mathcal{A}^+$ for every $\intervala.$ 
\end{enumerate}
\end{assu}

For every $\theta$ in $\mathcal{A}^+$ we denote with $Q^\theta$ the probability measure  defined by:    $$\frac{dQ^\theta}{dP}=\frac{U'(X^\theta_1)X^\theta_1}{\mathbb{E}\pqa U'(X^\theta_1)X^\theta_1\pqc}.$$
 
 The utility maximization problem consists in finding a process $\pi$ in $\mathcal{A}^+$ maximizing expected utility from terminal wealth, i.e.:
\beqn \label{21}
\pi=\arg \max_{\theta \in \mathcal{A}^+}\mathbb{E}\pqa U(X^\theta_1)\pqc.
\eeqn

Problem (\ref{21}) is not trivial because of the uncertain nature of the processes $A$ and $V$ and the non zero quadratic variation of $A.$ Indeed, let us suppose that $\pqa A\pqc=0$ and that both $A$ and $V$ are deterministic. Then,  
 it is sufficient to consider  $$\sup_{\lambda \in \mathbb{R}} \mathbb{E}\pqa U(X^\lambda_1)\pqc=\lim_{x\rightarrow +\infty}U(x),$$ and remind that $U$ is strictly increasing,  to see that a maximum can not be realized. 
The problem is less clear when the  term $-\frac{1}{2}\int_0^\cdot \theta_t^2d\pqa A\pqc_t$  and a source of randomness are added.

In the sequel, we will always assume the following.
\begin{assu}\label{5.42}
For every $\theta$ in $\cal{A}^+,$ 
$$
\mathbb{E}\pqa \vaa \int_0^1 \theta_td^- (A_t-V_t)\vac +\frac{1}{2}\int_0^1 \theta_t^2 \pqa A\pqc_t  \pqc<+\infty.
$$
\end{assu}

\begin{defi}
A process $\pi$ is said \textit{\textbf{optimal portfolio}} in $\cal{A}^+,$ if it is optimal for the function $\theta \mapsto \mathbb{E}\pqa U(X^\theta_1)\pqc$ in $\mathcal{A}^+,$ according to definition \ref{d3.18}.
\end{defi}

\begin{rema} Set $\xi=A-V,$ $\mathcal{A}=\mathcal{A}^+,$ and  
$$
F(\omega,x)=U\pta X_0(\omega)e^{ x+V_1(\omega)}\ptc, \quad (\omega,x) \in \Omega \times \mathbb{R}.
$$
According to definitions of section  \ref{s3.3.2}, $\mathcal{A}$ satisfies assumption \ref{a0}, the function $F$ is measurable, almost surely in $C^1(\mathbb{R}),$ strictly increasing and with bounded  first derivative. If $U$ satisfies assumption \ref{a5.41} then $F$ is also concave. Moreover $F(L^\theta)=U(X^\theta_1)$ for every $\theta$ in $\mathcal{A}^+.$  
\end{rema}

Before stating some results about the existence of an optimal portfolio, we provide examples of sets of admissible strategies with positive wealth.

\subsubsection{Admissible strategies via It\^o fields} 
For this example the reader should keep in mind subsection \ref{s4.0}.
 
\begin{prop}
Let $\cal{A}^+$ be the set of all processes $(\theta_t, \intervala)$ such that $\theta$ is the restriction to $[0,1)$ of a process $h$ belonging to $\mathcal{S}(\mathcal{C}^1_A(\mathbb{G})).$ Then $\mathcal{A}^+$ satisfies the hypotheses of assumption \ref{a5.40}.
\end{prop}

\begin{proof}
Let $h$ be in $\mathcal{S}(\mathcal{C}_A^1(\mathbb{G}))$ and $\theta$ its restriction on $[0,1).$ It is clear that $\theta$ is in $C^-_b([0,1))$ and $\mathbb{G}$-adapted. Thanks to proposition \ref{3.8}, $h$ is $A$-forward integrable, 
$\int_0^\cdot h_t d^-A_t$ belongs to  $\mathcal{S}(\mathcal{C}_A^2(\mathbb{G}))$ 
and the process $\int_0^\cdot h_t d^-A_t$ has finite quadratic variation equal to $\int_0^\cdot h_t^2 d\pqa A \pqc_t$. By remark \ref{r2.5}, $\int_0^\cdot h_td^-A_t=\int_0^\cdot \theta_t d^-A_t,$ and  conditions 1. and 2. of lemma \ref{l4.3} are thus satisfied. 
Remark \ref{r3.8} implies that the process  $$\exp \pta \int_0^\cdot \theta_td^-A_t+\int_{0}^\cdot (1-\theta_s)dV_s-\frac{1}{2}\int_0^\cdot \theta_t^2 d\pqa A \pqc_t \ptc $$ belongs to  $\mathcal{S}(\mathcal{C}^1_A(\mathbb{G})).$ Then,  by proposition \ref{3.8}, again, $\theta$ fulfills also condition 3. of lemma \ref{l4.3}. By construction, $\theta I_{[0,t]}$ is an element of $\mathcal{A}^+$ for every $\intervala$ and this concludes the proof. 
\end{proof}
 
\subsubsection{Admissible strategies via Malliavin calculus} \label{s5.4.2}
We restrict ourselves to the setting of section \ref{5.25}. We recall that in that case $A=\int_0^\cdot \sigma_tdW_t+\int_0^\cdot \mu_t dt.$ We make the following additional assumption:

$$S^0=e^{\int_0^\cdot r_tdt},
$$ 
with $r$ in $L^{1,z}$ for some $z>4$ and $\mathbb{F}$-adapted.
 
\begin{prop} \label{p5.29}
Let $\cal{A}^+$ be the set of all $\mathbb{G}$-adapted processes in  $C^-_b([0,1))$ being the restriction on $[0,1)$ of processes $h$ in $L^{1,2}_{-}\cap L^{2,2},$  such that $D^-h$ is in $L^{1,2}_{-},$ and the random variable  
$$
\sup_{t\in [0,1]} \pta \vaa h_t \vac+\sup_{s\in[0,1]} \vaa D_s h_t\vac+\sup_{s,u\in[0,1]} \vaa D_s D_u h_t\vac \ptc 
$$
is bounded. 

Then $\cal{A}^+$ satisfies the hypotheses of assumption \ref{a5.40}. 
\end{prop}

\begin{proof}
Let $h$ be as in the hypotheses. Proposition \ref{p4.2a} applies to to state  that $h$ is $A$-forward integrable and
\beqn \nonu
\int_0^\cdot h_td^-A_t&=&\int_0^\cdot h_t \sigma_t d^-W_t+\int_0^\cdot h_t\mu_t dt \\ \nonu
&=&\int_0^\cdot h_t \sigma_t \delta W_t+\int_0^\cdot \pta h_t\mu_t+ \sigma_tD_t^-h_t \ptc dt.
\eeqn
On the other hand, proposition \ref{c4.1} applies to obtain  
$$
\pqa  \int_0^\cdot h_td^-A_t \pqc= \pqa \int_0^\cdot h_t \sigma_t \delta W_t \pqc= \int_0^\cdot h_t^2\sigma_t^2dt.
$$
In particular, if $\theta$ is the restriction of $h$ on $[0,1),$ then $\theta$ fulfills point 1. and 2. of lemma \ref{l4.3}.

Consider the vector process $\pta \int_0^\cdot h_t d^-A_t, \int_0^\cdot (1-h_t)dV_t,\int_0^\cdot h_t^2 d\pqa A \pqc\ptc_t.$ We affirm that each of its components belongs to $L^{1,v}_-$ for some $v>4.$  In fact, the first component is equal to the sum of $\int_0^\cdot h_t\sigma_t\delta W_t$ and $\int_0^\cdot\pta \sigma_tD_t^-h_t + h_t\mu_t \ptc dt;$ the first term of the sum belongs to $L^{1,p}_-$ by lemma \ref{l4.6}, which applies thanks to remark \ref{r4.1}; the second term is in $L^{1,q \wedge p}_-$ thanks to lemma \ref{l4.4}; remark \ref{r4.1} and lemma \ref{l4.4} again imply that both $\int_0^\cdot (1-h_t)r_tdt$ and $\int_0^\cdot h_t^2 \sigma_t^2 dt$ belong to $L^{1,z}_-.$ 
We can thus apply proposition \ref{p4.15} to find that
$$
\int_0^\cdot X^h_t d^- \int_0^t h_s\sigma_sdW_s=\int_0^\cdot X^h_t h_t d^- \int_0^t \sigma_sdW_s,
$$
with $X^h=\exp \pta  \int_0^\cdot h_td^-A_t-\int_{0}^\cdot (1-h_s)dV_s-\frac{1}{2}\int_0^\cdot h^2 d\pqa A \pqc_t \ptc.$ 
This permits to conclude the proof. 
\end{proof}
 
\subsubsection{Admissible strategies via substitution}
We return here to the framework of subsection \ref{e5.21}. 
\begin{prop}
Let $\mathcal{A}^+$ be the set of all processes which are the restriction to $[0,1)$ of processes in $\mathcal{S}(\mathcal{A}^{p,\gamma}(L))$ for some $p>3$ and $\gamma>3d.$ Then $\mathcal{A}^+$ satisfies the hypotheses of assumption \ref{a5.40}.  
\end{prop}
\begin{proof}
Let $h$ be in $\mathcal{S}(\mathcal{A}^{p,\gamma}(L))$ for some $p>3$ and $\gamma >3d.$ Proposition \ref{p3.35} insures that $h$ is $A$-forward integrable, and that $\int_0^\cdot h _td^-A_t$  has finite quadratic variation equal to 
$
\int_0^\cdot h^2_td \pqa A\pqc_t.
$
The process $$X^h=\exp \pta \int_0^\cdot h_td^-A_t-\int_{0}^\cdot (1-h_t)dV_t-\frac{1}{2}\int_0^\cdot h^2_t d\pqa A \pqc_t\ptc$$ has bounded paths. Then, thanks to point 1. of remark \ref{r4.9}, to prove that 
\beqn \label{e}
\int_0^\cdot X^h_td^-\int_0^t h_s d^-A_s=\int_0^\cdot X^h_t h_t d^-A_t,
\eeqn
we are allowed to replace $X^h$ by $\psi(\log(X^h)),$ being $\psi$ a function of class $C^\infty(\mathbb{R})$ with bounded derivative.
Using lemma \ref{l4.17} it is possible to show that the process  $\psi \pta \log(X^h)\ptc$ belongs to $\mathcal{S}(\mathcal{A}^{\frac{p}{2},\frac{\gamma}{2}}(L))$. Proposition \ref{p3.35} again let us get equality (\ref{e}).
\end{proof}
 
\subsubsection{Optimal portfolios and $\mathcal{A}^+$-martingale property}
Adapting results contained in section \ref{s3.3.2} to the utility maximization problem, we can formulate the following propositions. We omit their proofs, being particular cases of the ones contained in that section.
  
\begin{prop}\label{p5.45}
If a process $\pi$ in $\mathcal{A}^+$ is an optimal portfolio, then  the process 
$
A-V-\int_0^\cdot \pi_td\pqa A\pqc_t
$
is an $\mathcal{A}^+$-martingale under ${Q^\pi}.$  If $U$ fulfills  assumption \ref{a5.41}, then the converse holds.
\end{prop}

\begin{prop} \label{p5.46} 
Suppose that $\mathcal{A}^+$ satisfies assumption $\mathcal{D}$ with respect to $\mathbb{G}$. 
If a process $\pi$ in $\cal{A}^+$ is an optimal portfolio, then the process 
$A-V-\int_0^\cdot \beta_t d\pqa A \pqc_t
$
is a $\mathbb{G}$-martingale under $P,$ with   $$\beta=\pi+\frac{1}{p^\pi}
\frac{d\pqa p^\pi,A\pqc}{d\pqa A \pqc},\quad \mbox{ and   }
\quad p^\pi=\mathbb{E}^{Q^\pi}\pqa \frac{dP}{dQ^\pi}\left|\right. \mathcal{G}_\cdot\pqc.$$
If $U$ fulfills  assumption \ref{a5.41}, then the converse holds.
\end{prop}

\begin{rema}
\begin{enumerate}
\item
We emphasize that if $U(x)=\log(x),$ then
the probability measure $Q^\pi$ appearing in propositions \ref{p5.45} and \ref{p5.46} is equal to $P.$
\item
 In \cite{AI} it is proved that if the maximum of expected logarithmic utility over all  \textit{simple admissible strategies} is finite, then $S$ is a semimartingale with respect $\mathbb{G}.$ This result does not imply proposition \ref{p5.46}. Indeed, we do not need to assume that our set of portfolio strategies contains the set of simple predictable admissible ones. On the contrary, we want to point out that, as soon as the class of admissible strategies is not \textit{large} enough, the semimartingale property of price processes could fail, even under finite expected utility.   
 \end{enumerate}
\end{rema}

\begin{prop}\label{p5.48}
Suppose that $U(x)=\log(x),$ $x$ in $(0,+\infty).$ Assume that there exists a measurable process $\gamma$ such that $A-V-\int_0^\cdot \gamma_t d\pqa A\pqc_t$ is an $\mathcal{A}^+$-martingale and there exists a sequence $(\theta^n)_{n\in \mathbb{N}}\subset \mathcal{A}^+$ such that
$$
\lim_{n \rightarrow +\infty}\mathbb{E}\pqa \int_0^1 \vaa \theta^n_t-\gamma_t\vac^2 d\pqa A\pqc_t \pqc=0.
$$ 
Then if $\gamma$ belongs to $\mathcal{A}^+$ it is an optimal portfolio. On the contrary, if an optimal portfolio $\pi$ exists, then 
$d\vaa \pqa A\pqc\vac \pga t\in [0,1),\pi_t \neq \gamma_t \pgc=0$  almost surely.  
\end{prop}

\subsubsection{Example}
We adopt the setting of section $\ref{s5.4.2}$ and we further assume that  $\sigma$ is a strictly positive real 

\begin{prop}
If a process $\pi$ is an optimal portfolio in $\cal{A}^+,$ then the process 
$
W-\int_0^\cdot \pta \frac{r_t-\mu_t}{\sigma}+\pi_t\sigma\ptc  dt
$
is an $\cal{A}^+$-martingale under $Q^\pi.$ If $U$ fulfills assumption \ref{a5.41}, then the converse holds.
\end{prop}
\begin{proof}
First of all we observe that it is not difficult to prove that $\mathcal{A}^+$ satisfies assumption \ref{5.42}.
If a process $\pi$ is an optimal portfolio in $\cal{A}^+$ then proposition \ref{p5.45} implies that the process $M^\pi$, with $M^\pi=\sigma\pta W-\int_0^\cdot \pta \frac{r_t-\mu_t}{\sigma}-\pi_t\sigma \ptc dt\ptc,$ is an $\cal{A}^+$-martingale under $Q^\pi.$
We observe that $\sigma^{-1}\mathcal{A}^+=\mathcal{A}^+.$ Therefore, $\sigma^{-1}M^\pi=W-\int_0^\cdot \pta \frac{r_t-\mu_t}{\sigma}+\pi_t\sigma\ptc  dt $ is an $\mathcal{A}^+$-martingale.

Similarly, if $U$ satisfies assumption \ref{a5.41}, the converse follows by proposition \ref{p5.45}. 
\end{proof}

\begin{coro}
Let $\mathcal{A}^+$ satisfy assumption $\mathcal{D}$ with respect to $\mathbb{G}.$ 
If a process $\pi$ in $\cal{A}^+$ is an optimal portfolio then the process 
$
\widetilde{W}=W-\int_0^\cdot \alpha_t  dt
$ with 
$$
\alpha=\pi\sigma+\frac{r-\mu}{\sigma}+\frac{1}{p^\pi}\frac{d\pqa p^\pi,W\pqc}{d\pqa W\pqc},\quad  \mbox{ and   }\quad p^\pi=\mathbb{E}^{Q^\pi}\pqa \frac{dP}{dQ^\pi}\left|\right. \mathcal{G}_\cdot\pqc,
$$
is a $\mathbb{G}$-Brownian motion under $P.$
If  $U$ satisfies assumption \ref{a5.41}, then the converse holds.
\end{coro}
\begin{proof}
Let $\pi$ be an optimal portfolio. By proposition \ref{p5.1}, the process $\widetilde{W}$ is a $\mathbb{G}$-martingale and so a $\mathbb{G}$-Brownian motion under $P.$
\end{proof}
The results concerning the example above were proved in \cite{BO}. We generalize them in two directions: we replace the geometric Brownian motion $A$ by a finite quadratic variation process and we let the set of possible strategies vary in sets which can, a priori, exclude some simple predictable processes.

\subsubsection{Example}
We consider the example treated in section \ref{3.26}. We suppose, for simplicity, that
$$
S_t=S_0 e^{ \sigma W_t+ \pta \mu-\frac{\sigma^2}{2}\ptc t}, \quad S^0_t=e^{rt}
\quad \interval,
$$
being $\sigma,$ $\mu$ and $r$ positive constants.  This implies $A_t=\sigma W_t+\mu t,$ and $V_t=rt$ for $\interval.$ 
We set $\mathcal{A}^+=\Theta(L).$

\begin{prop}
Suppose that $U(x)=\log(x),$ $x$ in $(0,+\infty).$ Suppose that $h(\cdot,L)$  belongs to the closure of $\Theta(L)$ in $L^2(\Omega \times [0,1]).$ Then an optimal portfolio $\pi$ exists if and only if the process $h(\cdot,L)+\int_0^\cdot \frac{\mu-r}{\sigma}dt $ belongs to $\Theta(L)$ and $\pi=h(\cdot,L)+\frac{\mu-r}{\sigma}.$ 
\end{prop}

\begin{proof}
The result follows from corollary \ref{c4.31}. 
\end{proof}

Sufficiency for the proposition above was shown, with more general $\sigma,$ $r$ and $\mu$ in theorem 3.2 of \cite{LNN}. Nevertheless, in this paper we go further in the analysis of utility maximization problem. Indeed, besides observing that the converse of that theorem holds true, we find that the existence of an optimal strategy is strictly connected, even for different choices of the utility function, to the $\cal{A}^+$-semimartingale property of $W.$ To be more precise, in that paper the authors show that an optimal process exists, under the given hypotheses, handling directly the expression of the expected utility, which has, in the logarithmic case, a nice expression. Here we reinterpret their techniques at a higher level which permits us to partially generalize those results.

\addcontentsline{toc}{chapter}{Bibliography}
\nocite{*}
\bibliographystyle{plain}
\bibliography{insider}

\end{document}